\documentclass{article}
\usepackage{latexsym}
\usepackage{amssymb}
\usepackage{amsmath}
\usepackage[english]{babel}
\usepackage{theorem}
\usepackage[mathscr]{eucal}
\usepackage{enumerate}

\input xypic
\xyoption{all}
\xyoption{arc}

\newtheorem{theo}{Theorem}[section]

\newtheorem{prop}[theo]{Proposition}

\newtheorem{lem}[theo]{Lemma}

\newtheorem{remark}[theo]{Remark}

\newcounter{equat}

\def\equat{\refstepcounter{equat}$$~}
\def\endequat{\leqno{\boldsymbol{(\arabic{equat})}}~$$}

\def\ov{\overline}

\def\alpt{{\tilde{\alpha}}}
\def\Srm{{\mathrm{S}}}
\def\Rrm{{\mathrm{R}}}
\def\Irm{{\mathrm{I}}}
\def\Cba{{\bar{C}}}

  \def\bG{{\mathfrak b}}  
    \def\CM{{\mathbb{C}}}

    \def\FM{{\mathbb{F}}}
  \def\gG{{\mathfrak g}}  \def\GM{{\mathbb{G}}}

  \def\lG{{\mathfrak l}}  
    
    \def\NM{{\mathbb{N}}}
    
    \def\PM{{\mathbb{P}}}
    \def\QM{{\mathbb{Q}}}
    \def\RM{{\mathbb{R}}}
    
  \def\tG{{\mathfrak t}}

    \def\ZM{{\mathbb{Z}}}

    \def\BC{{\mathcal{B}}}
    
    \def\DC{{\mathcal{D}}}
    \def\EC{{\mathcal{E}}}
    \def\FC{{\mathcal{F}}}

    \def\IC{{\mathcal{I}}}

    \def\LC{{\mathcal{L}}}
    
    \def\NC{{\mathcal{N}}}
    \def\OC{{\mathcal{O}}}

\def\Iti{{\tilde{I}}}

\def\a{\alpha}
\def\b{\beta}
\def\g{\gamma}

\def\D{\Delta}

\def\l{\lambda}

\def\s{\sigma}

\def\t{\tau}

\DeclareMathOperator{\Ad}{{\mathrm{Ad}}}
\DeclareMathOperator{\Coker}{{\mathrm{Coker}}}

\DeclareMathOperator{\Hom}{{\mathrm{Hom}}}
\DeclareMathOperator{\Ker}{{\mathrm{Ker}}}

\DeclareMathOperator{\Tr}{{\mathrm{Tr}}}

\newcommand{\elem}[1]{\stackrel{#1}{\longto}}

\def\iff{\Leftrightarrow}

\def\to{\rightarrow}
\def\longto{\longrightarrow}

\def\fonctio#1#2#3#4{\begin{array}{ccc}
{#1} & \longto & {#2} \\
{#3} & \longmapsto & {#4} 
\end{array}}

\def\incl{\hspace{0.05cm}{\subset}\hspace{0.05cm}}

\def\lexp#1#2{\kern\scriptspace\vphantom{#2}^{#1}\kern-\scriptspace#2}
\def\le{\hspace{0.1em}\mathop{\leqslant}\nolimits\hspace{0.1em}}
\def\ge{\hspace{0.1em}\mathop{\geqslant}\nolimits\hspace{0.1em}}

\mathchardef\inferieur="321E
\mathchardef\superieur="321F

\def\eqna{\begin{eqnarray*}}
\def\endeqna{\end{eqnarray*}}

\def\mini{{\mathrm{min}}}

\def\iff{if and only if}

\newcommand{\proof}[1]{\noindent {\bf Proof~:}\hskip.5cm #1~\(\Box\)\\}

\def\longue{{\mathrm{lg}}}
\def\courte{{\mathrm{sh}}}
\def\hauteur{{\mathrm{ht}}}

\author{Daniel Juteau\thanks{UFR de math\'ematiques, Universit\'e 
Denis Diderot Paris 7}}
\title{Cohomology of the minimal nilpotent orbit}

\begin{document}

\maketitle

\begin{abstract}
We compute the integral cohomology of the minimal non-trivial
nilpotent orbit in a complex simple (or quasi-simple) Lie algebra.
We find by a uniform approach that the middle cohomology group
is isomorphic to the fundamental group of the 
sub-root system generated by the long simple roots.
The modulo $\ell$ reduction of the Springer correspondent
representation involves the sign representation exactly
when $\ell$ divides the order of this cohomology group.
The primes dividing the torsion of the rest of the cohomology are
bad primes.
\end{abstract}

\section*{Introduction}

Let $G$ be a quasi-simple complex Lie group,
with Lie algebra $\gG$.
We denote by $\NC$ the nilpotent variety of $\gG$.
The group $G$ acts on $\NC$ by the adjoint action,
with finitely many orbits. If $\OC$ and $\OC'$ are two orbits,
we write $\OC \le \OC'$ if $\OC \incl \ov{\OC'}$. This defines
a partial order on the adjoint orbits. It is well known that
there is a unique minimal non-zero orbit $\OC_\mini$ (see for example
\cite{CM}, and the introduction of \cite{KP2}). 
The aim of this article is to compute the integral cohomology
of $\OC_\mini$.

The nilpotent variety $\NC$ is a cone in $\gG$: it is closed under
multiplication by a scalar. Let us consider its image $\PM(\NC)$
in $\PM(\gG)$. It is a closed subvariety of this projective space,
so it is a projective variety. Now $G$ acts on $\PM(\NC)$,
and the orbits are the $\PM(\OC)$, where $\OC$ is a non-trivial adjoint
orbit in $\NC$. The orbits of $G$ in $\PM(\NC)$ are ordered in the same
way as the non-trivial orbits in $\NC$. Thus $\PM(\OC_\mini)$ is the
minimal orbit in $\PM(\NC)$, and therefore it is closed:
we deduce that it is a projective variety.
Let $x_\mini \in \OC_\mini$, and let $P = N_G(\CM x_\mini)$
(the letter $N$ stands for normalizer, or setwise stabilizer).
Then $G/P$ can be identified to $\PM(\OC_\mini)$, which is a projective
variety. Thus $P$ is a parabolic subgroup of $G$.
Now we have a resolution of singularities (see section \ref{res spectral})
\[
G \times_P \CM x_\mini \longto \ov{\OC_\mini} = \OC_\mini \cup \{0\}
\]
which restricts to an isomorphism
\[
G \times_P \CM^* x_\mini \elem{\sim} \OC_\mini.
\]

From this isomorphism, one can already deduce that the dimension of
$\OC_\mini$ is equal to one plus the dimension of $G/P$. If we
fix a maximal torus $T$ in $G$ and a Borel subgroup $B$ containing it,
we can take for $x_\mini$ a highest weight vector for the adjoint action
on $\gG$. Then $P$ is the standard parabolic subgroup corresponding
to the simple roots orthogonal to the highest root, and the dimension
of $G/P$ is the number of positive roots not orthogonal to the highest
root, which is $2h - 3$ in the simply-laced types, where
$h$ is the Coxeter number (see \cite[chap. VI, \S 1.11, prop. 32]{BOUR456}).
So the dimension of $\OC_\mini$ is $2h - 2$ is that case. In \cite{WANG},
Wang shows that this formula is still valid if we replace $h$ 
by the dual Coxeter number $h^\vee$ (which is equal to $h$ only in the
simply-laced types).

We found a similar generalization of a result of Carter
(see \cite{CARTER}), relating
the height of a long root to the length of an element of minimal length
taking the highest root to that given long root, in the simply-laced
case: the result extends to all types, if we take the height of the
corresponding coroot instead (see Section
\ref{philg xi}, and Theorem \ref{conclu racines}). 

To compute the cohomology of $\OC_\mini$,
we will use the Gysin sequence
associated to the $\CM^*$-fibration $G \times_P \CM^* x_\mini \longto G/P$.
The Pieri formula of Schubert calculus gives an answer in terms
of the Bruhat order (see section \ref{res spectral}).
Thanks to the results of section \ref{philg xi},
we translate this in terms of the combinatorics of the root system
(see Theorem \ref{spectral roots}). As a consequence, we obtain the
following results (see Theorem \ref{middle}):

\medskip

\noindent{\bf Theorem~} 
{\it
(i)~
The middle cohomology of $\OC_\mini$ is given by
$$H^{2h^\vee - 2}(\OC_\mini, \ZM) \simeq P^\vee({\Phi'})/Q^\vee({\Phi'})$$
where $\Phi'$ is the sub-root system of $\Phi$ generated by the long
simple roots,
and $P^\vee({\Phi'})$ (resp. $Q^\vee({\Phi'})$) is its coweight lattice
(resp. its coroot lattice).

(ii)~
If $\ell$ is a good prime for $G$, then there is no $\ell$-torsion
in the rest of the cohomology of $\OC_\mini$.
}

\medskip

Part $(i)$ is obtained by a general argument, while $(ii)$ is obtained
by a case-by-case analysis (see section \ref{case by case},
where we give tables for each type).

In section \ref{a bis},
we explain a second method for the type $A_{n - 1}$, based on another
resolution of singularities: this time, it is a cotangent bundle
on a projective space (which is also a generalized flag variety).
This cannot be applied to other types, because the minimal class
is a Richardson class only in type $A$.

The motivation for this calculation is the modular
representation theory of the Weyl group $W$.
To each rational irreducible representation of $W$,
one can associate, via the Springer correspondence
(see for example \cite{SPRTRIG, SPRWEYL, BM, KL, SLO1, ICC, SHOJI}),
a pair consisting in a nilpotent
orbit and a $G$-equivariant local system on it (or, equivalently,
a pair $(x,\chi)$ where $x$ is a nilpotent element of $\gG$,
and $\chi$ is an irreducible character of the finite group
$A_G(x) = C_G(x) / C_G^0(x)$, up to $G$-conjugation).
Note that Springer's construction differs from the others by
the sign character.
All the pairs consisting of a nilpotent orbit and the constant
sheaf on this orbit arise in this way. In the simply-laced
types, the irreducible representation of $W$ corresponding
to the pair $(\OC_\mini, \QM)$ is the natural representation
tensored with the sign representation. In the other types,
we have a surjection from $W$ to the reflection group $W'$
corresponding to the subdiagram of the Dynkin diagram of $W$
consisting in the long simple roots. The Springer correspondent
representation is then the natural representation of $W'$
lifted to $W$, tensored with the sign representation.

We believe that the decomposition matrix of the Weyl group
(and, in fact, of an associated Schur algebra) can be deduced
from the decomposition matrix of $G$-equivariant perverse sheaves on
the nilpotent variety $\NC$.
In \cite{DEC}, we will use Theorem \ref{middle} to determine
some decomposition numbers for perverse sheaves (which give some evidence
for this conjecture).
Note that we are really interested in the torsion.
The rational cohomology must already be known to the experts
(see Remark \ref{rat}).

All the results and proofs of this article remain valid for
$G$ a quasi-simple reductive group over $\overline{\FM_p}$,
with $p$ good for $G$, using the \'etale topology. In this context, 
one has to take $\QM_\ell$ and $\ZM_\ell$ coefficients,
where $\ell$ is a prime different from $p$, instead
of $\QM$ and $\ZM$.

\section{Long roots and distinguished coset representatives}\label{philg xi}

The Weyl group $W$ of an irreducible and reduced root system $\Phi$
acts transitively on the set $\Phi_\longue$ of long roots in $\Phi$,
hence if $\a$ is an element of $\Phi_\longue$,
then the long roots are in bijection with $W/W_\a$,
where $W_\a$ is the stabilizer of $\a$ in $W$ (a parabolic subgroup).
Now, if we fix a basis $\D$ of $\Phi$, and if we choose for $\a$
the highest root $\alpt$, we find a relation between the partial orders
on $W$ and $\Phi_\longue$ defined by $\D$, and between the length
of a distinguished coset representative and the (dual) height
of the corresponding long root. After this section was written,
I realized that the result was already proved by Carter in the
simply-laced types in \cite{CARTER} (actually, this result is quoted
in \cite{SPRTRIG}). We extend it to any type and
study more precisely the order relations involved.
I also came across \cite[\S 4.6]{BB}, where the depth of a positive
root $\b$ is defined as the minimal integer $k$ such that
there is an element $w$ in $W$ of length $k$ such that $w(\b) < 0$.
By the results of this section, the depth of a positive long root is nothing
but the height of the corresponding coroot (and the depth
of a positive short root is equal to its height).

For the classical results about root systems that are used throughout
this section, the reader may refer to \cite[Chapter VI, \S 1]{BOUR456}.
It is now available in English \cite{BOUR456ENG}.

\subsection{Root systems}\label{root systems} 
Let $V$ be a finite dimensional $\RM$-vector space and
$\Phi$ a root system in $V$.
We note $V^*=\Hom(V,\RM)$ and, if $\a \in \Phi$, we denote by $\a^\vee$ 
the corresponding coroot
and by $s_\a$ the reflexion $s_{\a,\a^\vee}$ of 
\cite[chap. VI, \S 1.1, d\'ef. 1, $(\Srm\Rrm_{\Irm\Irm})$]{BOUR456}. 
Let $W$ be the Weyl group of $\Phi$. 
The perfect pairing between $V$ and $V^*$ will be denoted by $\langle,\rangle$. 
Let $\Phi^\vee=\{\a^\vee~|~\a \in \Phi\}$. 
In all this section, we will assume that $\Phi$ is {\it irreducible} 
and {\it reduced}. 
Let us fix a scalar product $(~|~)$ on $V$, invariant under $W$,
such that
$$\min_{\a \in \Phi} (\a|\a) = 1.$$
We then define the integer 
\[
r=\max_{\a \in \Phi} (\a|\a).
\]
Let us recall that, since $\Phi$ is irreducible and reduced,
we have $r \in \{1,2,3\}$ 
and $(\a|\a) \in \{1,r\}$ if $\a \in \Phi$ 
(see \cite[chap. VI, \S 1.4, prop. 12]{BOUR456}). We define
\[
\Phi_\longue=\{\a \in \Phi~|~(\a|\a) = r\}
\]
and
\[
\Phi_\courte=\{\a \in \Phi~|~(\a|\a) < r \}
= \Phi\setminus \Phi_\longue.
\]
If $\a$ and $\b$ are two roots, then 
\equat\label{scalaire}
\langle \a,\b^\vee \rangle = \frac{2(\a|\b)}{(\b|\b)}.
\endequat
In particular, if $\a$ and $\b$ belong to $\Phi$, then
\equat\label{entier}
2(\a|\b) \in \ZM
\endequat
and, if $\a$ or $\b$ belongs to $\Phi_\longue$, then
\equat\label{r entier}
2(\a|\b) \in r\ZM
\endequat
The following classical result says that $\Phi_\longue$ is a closed
subset of $\Phi$.

\begin{lem}\label{somme longue}
If $\a$, $\b \in \Phi_\longue$ are such that $\a+\b \in \Phi$, 
then $\a+\b \in \Phi_\longue$.
\end{lem}

\proof{ 
We have $(\a+\b~|~\a+\b)=(\a|\a)+(\b|\b)+2(\a|\b)$.
Thus, by \ref{r entier}, 
we have $(\a+\b~|~\a+\b) \in r\ZM$, which implies the desired result.
}

\subsection{Basis, positive roots, height}\label{basis}
Let us fix a basis $\D$ of $\Phi$ and let $\Phi^+$ be the set of roots
$\a \in \Phi$ whose coefficients in the basis $\D$ are non-negative. 
Let $\D_\longue = \Phi_\longue \cap \D$ and
$\D_\courte = \Phi_\courte \cap \D$. 
Note that $\D_\longue$ need not be a basis of $\Phi_\longue$. 
Indeed, $\Phi_\longue $ is a root system of rank equal to the rank
of $\Phi$, whereas $\D_\longue$ has fewer elements than $\D$ if
$\Phi$ is of non-simply-laced type.
Let us recall the following well-known result
\cite[chap. VI, \S 1, exercice 20 (a)]{BOUR456}:

\begin{lem}\label{caracterisation longue}
Let $\g \in \Phi$ and write $\g=\sum_{\a \in \D} n_\a \a$, with
$n_\a \in \ZM$. Then $\g \in \Phi_\longue$ \iff \
$r$ divides all the $n_\a$, $\a \in \D_\courte$.
\end{lem}

\proof{
Let $\Phi'$ be the set of roots $\g' \in \Phi$ such that,
if $\g'=\sum_{\a \in \D} n_\a' \a$, then $r$ divides $n_\a'$ for all
$\a \in \D_\courte$. We want to show that $\Phi_\longue = \Phi'$. 

Suppose that $r$ divides all the $n_\a$, $\a \in \D_\courte$. Then
$n_\a^2 (\a|\a) \in r\ZM$ for all $\a\in\D$, and by
\ref{entier} and \ref{r entier}, we have $2n_\a n_\b(\a|\b) \in r\ZM$ 
for all $(\a,\b) \in \D \times \D$ such that $\a \neq \b$. 
Thus $(\g|\g) \in r\ZM$, which implies that $\g \in \Phi_\longue$. 
Thus $\Phi' \subset \Phi_\longue$.

Since $W$ acts transitively on $\Phi_\longue$, 
it suffices to show that $W$ 
stabilizes $\Phi'$. In other words, it is enough to show that,
if $\a \in \D$ and $\g \in \Phi'$, then $s_\a(\g) \in \Phi'$. But
$s_\a(\g)=\g -\langle \g,\a^\vee \rangle \a$. 
If $\a \in \D_\longue$, then $s_\a(\g) \in \Phi'$ because $\g \in \Phi'$. 
If $\a \in \D_\courte$, then $\langle \g,\a^\vee\rangle = 2(\g|\a) \in r\ZM$ 
because $\g \in \Phi' \subset \Phi_\longue$
(see \ref{scalaire} and \ref{r entier}). 
Thus $s_\a(\g) \in \Phi'$.
}

If $\g = \sum_{\a \in \D} n_\a \a \in \Phi$, the {\it height} of $\g$ 
(denoted by $\hauteur(\g)$) is defined by
$\hauteur(\g)=\sum_{\a \in \D} n_\a$. 
One defines the height of a coroot similarly.

If $\g$ is long, we have
$$\g^\vee =
\sum_{\a \in \D_\longue} n_\a \a^\vee
 + \frac{1}{r} \sum_{\a\in\D_\courte} n_\a \a^\vee.$$
Let
$$\hauteur^\vee(\g) := \hauteur(\g^\vee)
= \sum_{\a \in \D_\longue} n_\a + \frac{1}{r} \sum_{\a\in\D_\courte} n_\a.$$
In particular, the right-hand side of the last equation is an integer,
which is also a consequence of Lemma \ref{caracterisation longue}.

If $\a$ and $\b$ are long roots such that $\a + \b$ is a
(long) root, then $(\a + \b)^\vee = \a^\vee + \b^\vee$, so
$\hauteur^\vee$ is additive on long roots.

\subsection{Length}\label{length}
Let $l : W \to \NM=\{0,1,2,\dots\}$ be the {\it length} function associated 
to $\D$: if we let
$$N(w)=\{\a \in \Phi^+~|~w(\a) \in - \Phi^+\},$$
then we have 
\equat\label{longueur}
l(w)=|N(w)|.
\endequat
If $\a \in \Phi^+$ and if $w \in W$, then we have
\equat\label{ajout}
\text{\it $l(ws_\a) > l(w)$ \iff\  $w(\a) \in \Phi^+$.}
\endequat
Replacing $w$ by $w^{-1}$, and using the fact that an element
of $W$ has the same length as its inverse, we get
\equat\label{ajout gauche}
\text{\it $\ell(s_\a w) > \ell(w)$ \iff\  $w^{-1}(\a) \in \Phi^+$.}
\endequat
More generally, it is easy to show that,
if $x$ and $y$ belong to $W$, then
\equat\label{dyer}
N(xy) = N(y) \stackrel{.}{+} {}^{y^{-1}}N(x)
\endequat
where $\stackrel{.}{+}$ denotes the symmetric difference
(there are four cases to consider), and therefore
\equat\label{ajout bis}
\text{\it $l(xy)=l(x)+l(y)$ \iff\  $N(y) \subset N(xy)$}.
\endequat

Let $w_0$ be the longest element of $W$.
Recall that
\equat\label{wo}
l(w_0w)=l(ww_0)=l(w_0)-l(w)
\endequat
for all $w \in W$. If $I$ is a subset of $\D$, we denote by 
$\Phi_I$ the set of the roots $\a$ which belong
to the sub-vector space of $V$ generated by $I$ and we let
$$\Phi_I^+ = \Phi_I \cap\Phi^+ \qquad\text{and}\qquad
W_I=<s_\a~|~\a \in I>.$$
We also define
$$X_I=\{w \in W~|~w(\Phi_I^+) \subset \Phi^+\}.$$
Let us recall that $X_I$ is a set of coset representatives of $W/W_I$ and that
$w \in X_I$ \iff\  $w$ is of minimal length in $wW_I$. 
Moreover, we have
\equat\label{somme longueurs}
l(xw)=l(x)+l(w)
\endequat
if $x \in X_I$ and $w \in W_I$. 
We denote by $w_I$ the longest element of $W_I$. Then $w_0w_I$ is
the longest element of $X_I$ (this can be easily deduced from \ref{wo} 
and \ref{somme longueurs}). 
Finally, if $i$ is an integer, we denote by $W^i$ the set of elements
of $W$ of length $i$, and similarly $X_I^i$ is the set of elements
of $X_I$ of length $i$. To conclude this section,
we shall prove the following result, which should be well known:

\begin{lem}\label{longueur reflexion}
If $\b \in \Phi_\longue^+$, then
$l(s_\b) = 2~\hauteur^\vee(\b) - 1$. 
\end{lem}

\proof{We shall prove the result by induction on
$\hauteur^\vee(\b)$. The case where $\hauteur^\vee(\b) = 1$ is clear. 
Suppose $\hauteur^\vee(\b) > 1$ and suppose the result holds
for all positive long roots whose dual height is strictly smaller.

First, there exists a $\g \in \D$ such that $\b - \g \in \Phi^+$ 
(see \cite[chap. VI, \S 1.6, prop. 19]{BOUR456}). Let $\a = \b - \g$. 
There are two possibilities:

\medskip

$\bullet$ If $\g \in \D_\longue$, then $\a = \b - \g \in \Phi_\longue$
by Lemma \ref{somme longue}. 
Moreover, $\hauteur^\vee(\a) = \hauteur^\vee(\b) - 1$. Thus
$l(s_\a) = 2~\hauteur^\vee(\a)-1$. We have $(\a|\g) \neq 0$
(otherwise $\b = \a + \g$ 
would be of squared length $2r$, which is impossible). By 
\cite[chap. VI, \S 1.3]{BOUR456}, we have $(\a|\g) = - r/2$. 
Thus $\b = s_\g(\a) = s_\a(\g)$, and $s_\b = s_\g s_\a s_\g$. Since 
$s_\a(\g) > 0$, we have $l(s_\b s_\a) = l(s_\b) +1$ (see \ref{ajout}).
Since $s_\g s_\a (\g) = s_\g(\b) = \a > 0$, we have
$l(s_\g s_\a s_\g) = l(s_\a s_\g) + 1 = l(s_\a) + 2$ (see \ref{ajout}),
as expected.

\medskip

$\bullet$ If $\g \in \D_\courte$, then, by 
\cite[chap. VI, \S 1.3]{BOUR456}, we have 
$\a = \b - r\g \in \Phi_\longue^+$, $(\a|\g) = -r/2$, and 
$\hauteur^\vee(\a) = \hauteur^\vee(\b) - 1$.
As in the first case, we have $\b = s_\g(\a)$. Thus 
$s_\b = s_\g s_\a s_\g$ and the same argument applies.
}

\begin{remark}\label{short}
\emph{By duality, if $\b \in \Phi_\courte^+$, we have}
$$l(s_\b) = 2~\hauteur(\b) - 1.$$
\end{remark}

\subsection{Highest root}\label{highest root} 
Let $\alpt$ be {\it the highest root} of $\Phi$ relatively to $\D$ 
(see \cite[chap. VI, \S 1.8, prop. 25]{BOUR456}). 
It is of height $h - 1$, where $h$ is the Coxeter number of $\Phi$.
The \emph{dual Coxeter number} $h^\vee$ can be defined as
$1 + \hauteur^\vee(\alpt)$. Let us recall the following facts:
\equat\label{grand long}
\alpt \in \Phi_\longue
\endequat
and
\equat\label{positif}
\text{\it If $\a \in \Phi^+ \setminus\{\alpt\}$, 
then $\langle \a,\alpt^\vee \rangle \in \{0,1\}$.}
\endequat
In particular,
\equat\label{plus}
\text{\it If $\a \in \Phi^+$, then $\langle \alpt,\a^\vee \rangle \ge 0$}
\endequat
and
\equat\label{chambre grand}
\alpt \in \Cba,
\endequat
where $C$ is the chamber associated to $\D$.

From now on, ${\Iti}$ will denote the subset of $\D$ defined by
\equat\label{Iti}
{\Iti}=\{\a \in \D~|~(\alpt|\a)=0\}.
\endequat
By construction, ${\Iti}$ is stable under any automorphism of $V$
stabilizing $\D$. 
In particular, it is stable under $-w_0$. 
By \ref{plus}, we have
\equat\label{phii}
\Phi_{\Iti} = \{\a \in \Phi~|~(\alpt|\a) = 0\}.
\endequat
From \ref{phii} and \cite[chap. V, \S 3.3, prop. 2]{BOUR456}, we
deduce that
\equat\label{stab}
W_{\Iti}=\{w \in W~|~w(\alpt)=\alpt\}.
\endequat
Note that $w_0$ and $w_{\Iti}$ commute (because $-w_0({\Iti})= {\Iti}$). 
We have
\equat\label{wo wi}
N(w_0w_{\Iti}) = \Phi^+ \setminus \Phi_{\Iti}^+.
\endequat

Let us now consider the map $W \to \Phi_\longue$, $w \mapsto w(\alpt)$. 
It is surjective \cite[chap. VI, \S 1.3, prop. 11]{BOUR456} 
and thus induces a bijection $W/W_I \to \Phi_\longue$ by \ref{stab}. 
It follows that the map
\equat\label{bijection}
\fonctio{X_{\Iti}}{\Phi_\longue}{x}{x(\alpt)}
\endequat
is a bijection. If $\a \in \Phi_\longue$, we will denote by $x_\a$ the unique 
element of $X_{\Iti}$ such that $x_\a(\alpt)=\a$. We have
\equat\label{conjugaison}
x_\a s_\alpt = s_\a x_\a.
\endequat

\begin{lem}\label{salpt}
We have $w_0 w_{\Iti} = w_{\Iti} w_0 = s_\alpt$.
\end{lem}

\proof{
We have already noticed that $w_0$ and $w_{\Iti}$ commute.

In view of \cite[chap. VI, \S 1, exercice 16]{BOUR456}, 
it suffices to show that $N(w_0w_{\Iti}) = N(s_\alpt)$, 
that is, $N(s_\alpt) = \Phi^+ \setminus \Phi_{\Iti}^+$ (see \ref{wo wi}). 
First, if $\a \in \Phi_{\Iti}^+$, then $s_\alpt(\a) = \a$, so that 
$\a \not\in N(s_\alpt)$. This shows that $N(s_\alpt) \subset
\Phi^+ \setminus \Phi_{\Iti}^+$. 

Let us show the other inclusion. If $\a \in \Phi^+ \setminus \Phi_{\Iti}^+$,
then $\langle \alpt,\a^\vee \rangle > 0$ 
by \ref{plus} and \ref{phii}.
In particular, $s_\alpt(\a) = \a-\langle \a,\alpt^\vee \rangle \alpt$
cannot belong to $\Phi^+$ since $\alpt$ is the highest root.
}

\begin{prop}\label{palindrome}
Let $\a\in\Phi_\longue^+$. Then we have
$$l(x_\a s_\alpt) = l(s_\alpt) - l(x_\a)$$
\end{prop}

\proof{
We have
$$\begin{array}{rcll}
l(x_\a s_\alpt) &=& l(x_\a w_{\Iti} w_0) & \textrm{by Lemma \ref{salpt}}\\
&=& l(w_0) - l(x_\a w_{\Iti}) & \textrm{by \ref{wo}}\\
&=& l(w_0) - l(w_{\Iti}) - l(x_\a) & \textrm{by \ref{somme longueurs}}\\
&=& l(w_0w_{\Iti}) - l(x_\a) & \textrm{by \ref{wo}}\\
&=& l(s_\alpt) - l(x_\a) & \textrm{by Lemma \ref{salpt}.}
\end{array}$$
}

\begin{prop}\label{saxa}
If $\a\in\Phi_\longue^+$, then $x_{-\a} = s_\a x_\a$ and
$l(x_{-\a}) = l(s_\a x_\a) = l(s_\a) + l(x_\a)$.
\end{prop}

\proof{
We have $s_\a x_\a(\alpt) = s_\a(\a) = -\a$, so to show
that $x_{-\a} = s_\a x_\a$, it is enough to show that
$s_\a x_\a \in X_{\Iti}$. But, if $\b \in \Phi_{\Iti}^+$, we have
(see \ref{conjugaison}) $s_\a x_\a(\b) = x_\a s_\alpt(\b)
= x_\a(\b) \in \Phi^+$. Hence the first result.

Let us now show that the lengths add up. By \ref{ajout bis},
it is enough to show that $N(x_\a) \incl N(s_\a x_\a)$.
Let then $\b\in N(x_\a)$. Since $x_\a \in X_{\Iti}$, $\b$ cannot
be in $\Phi_{\Iti}^+$. Thus $\langle \b, \alpt^\vee \rangle > 0$.
Therefore, $\langle x_\a(\b), \a^\vee \rangle > 0$.
Now, we have
$s_\a x_\a(\b) = x_\a(\b) - \langle x_\a(\b), \a^\vee \rangle \a < 0$
(remember that $x_\a(\b) < 0$ since $\b\in N(x_\a)$). 
}

\begin{prop}\label{xa}
For $\a\in\Phi_\longue^+$, we have
$$l(x_\a) = \frac{l(s_\alpt) - l(s_\a)}{2}
          = \hauteur^\vee(\alpt) - \hauteur^\vee(\a)$$
$$l(x_{-\a}) = \frac{l(s_\alpt) + l(s_\a)}{2}
           = \hauteur^\vee(\alpt) + \hauteur^\vee(\a) - 1$$
\end{prop}

\proof{
This follows from Propositions \ref{palindrome} and \ref{saxa},
and \ref{conjugaison}.
}

\subsection{Orders}\label{orders}

The choice of $\D$ determines an order relation on $V$.
For $x$, $y\in V$, we have
$y \le x$ \iff\  $y - x$ is a linear combination
of the simple roots with non-negative coefficients.

For $\a\in\Phi_\longue$, it will be convenient to define the \emph{level}
$L(\a)$ of $\a$ as follows:
\equat\label{level}
L(\a) =
\left\{
\begin{array}{ll}
\hauteur^\vee(\alpt) - \hauteur^\vee(\a)
& \textrm{if $\a > 0$}\\
\hauteur^\vee(\alpt) - \hauteur^\vee(\a) - 1
& \textrm{if $\a < 0$}
\end{array}
\right.
\endequat
If $i$ is an integer, let $\Phi_\longue^i$ be the set of long
roots of level $i$. Then Proposition \ref{xa} says that
the bijection \ref{bijection} maps $X_\Iti^i$ onto $\Phi_\longue^i$.

For $\g\in\Phi^+$, we write
\equat\label{elem}
\b \elem{\g} \a \textrm{ \iff\  }
\a = s_\g(\b) \textrm{ and } L(\a) = L(\b) + 1.
\endequat
In that case, we have $\b - \a = \langle b, \g^\vee \rangle \g > 0$, so
$\b > \a$.

If $\a$ and $\b$ are two long roots, we say that
there is a \emph{path} from $\b$ to $\a$, and we write
$\a \preceq \b$, \iff\  there exists a sequence $(\b_0, \b_1, \ldots, \b_k)$ of
long roots,  and a sequence $(\g_1,\ldots,\g_k)$ of positive roots,
such that
\equat\label{path}
\b = \b_0 \elem{\g_1} \b_1 \elem{\g_2} \ldots \elem{\g_k} \b_k = \a.
\endequat
In that case, we have $L(\b_i) = L(\b) + i$ for 
$i \in \{0,~\ldots,~k\}$. If moreover all the roots $\g_i$ are simple,
we say that there is a \emph{simple path} from $\b$ to $\a$.

On the other hand, we have the Bruhat order on $W$ defined by
the set of simple reflections $S = \{s_\a~\mid~\a\in\D\}$.
If $w$  and $w'$ belong to $W$, we write
$w\longto w'$ if $w' = s_\g w$ and $l(w') = l(w) + 1$, for some
positive root $\g$. In that case, we write $w\elem{\g}w'$
(the positive root $\g$ is uniquely determined). The Bruhat
order $\le$ is the reflexive and transitive closure of the
relation $\longto$. On $X_\Iti$, we will consider the restriction of
the Bruhat order on $W$.

Let us now consider the action of a simple reflection on a long root.

\begin{lem}\label{action simple}
Let $\b\in\Phi_\longue$ and $\g\in\D$. Let $\a = s_\g(\b)$.

\begin{enumerate}
\item
\begin{enumerate}[(i)]
\item If $\b\in\D_\longue$ and $(\b|\g) > 0$, then $\g = \b$,
$\a = -\b$ and $\langle \b, \g^\vee \rangle = 2$.
\item If $\b\in -\D_\longue$ and $(\b|\g) < 0$, then $\g = -\b$,
$\a = -\b$ and $\langle \b, \g^\vee \rangle = - 2$.
\item Otherwise, $\a$ and $\langle \b, \g^\vee \rangle$ are
given by the following table:
$$
\begin{array}{|c|c|c|}
\hline
& \g\in\D_\longue & \g\in\D_\courte\\
\hline
\begin{array}{c}
(\b|\g) > 0\\
(\b|\g) = 0\\
(\b|\g) < 0
\end{array}
&
\begin{array}{l|l}
\a = \b - \g & \langle \b, \g^\vee \rangle = 1\\
\a = \b      & \langle \b, \g^\vee \rangle = 0\\
\a = \b + \g & \langle \b, \g^\vee \rangle = - 1\\
\end{array}
&
\begin{array}{l|l}
\a = \b - r\g & \langle \b, \g^\vee \rangle = r\\
\a = \b       & \langle \b, \g^\vee \rangle = 0\\
\a = \b + r\g & \langle \b, \g^\vee \rangle = - r\\
\end{array}\\
\hline
\end{array}
$$
\end{enumerate}

\item
\begin{enumerate}[(i)]
\item If $(\b|\g) > 0$ then $L(\a) = L(\b) + 1$, so that $\b \elem{\g} \a$.
\item If $(\b|\g) = 0$ then $L(\a) = L(\b)$, and in fact $\a = \b$.
\item If $(\b|\g) < 0$ then $L(\a) = L(\b) - 1$, so that $\a \elem{\g} \b$.
\end{enumerate}
\end{enumerate}
\end{lem}
 
\proof{
Part 1 follows from inspection of the possible cases in
\cite[Chapitre VI, \S 1.3]{BOUR456}.

Part 2 is a consequence of part 1. Note that there is a special case
when we go from positive roots to negative roots, and \emph{vice versa}.
This is the reason why there are two cases in the definition of the
level.
}

To go from a long simple root to the opposite of a long simple root,
one sometimes needs a non-simple reflection.

\begin{lem}\label{through}
Let $\b\in\D_\longue$, $\a\in -\D_\longue$ and $\g\in\Phi^+$.
Then $\b\elem{\g}\a$ \iff\  we are in one of the following cases:
\begin{enumerate}[(i)]
\item $\a = -\b$ and $\g = \b$. In this case,
$\langle \b, \g^\vee \rangle = 2$.
\item  $\b + (-\a)$ is a root and $\g = \b + (-\a)$. In this case,
$\langle \b, \g^\vee \rangle = 1$.
\end{enumerate}
\end{lem}

\proof{
This is straightforward.
}

But otherwise, one can use simple roots at each step.

\begin{prop}\label{simple path}
Let $\a$ and $\b$ be two long roots such that $\a \le \b$.
Write $\b = \sum_{\s\in J} n_\s \s$, and $\a = \sum_{\t\in K} m_\t \t$,
where $J$ (resp. $K$) is a non-empty subset of $\D$, and the $n_\s$
(resp. the $m_\t$) are non-zero integers, all of the same sign.
\begin{enumerate}[(i)]
\item\label{path positive} If $0 < \a \le \b$,
then there is a simple path from $\b$ to $\a$.
\item\label{path through} If $\a < 0 < \b$,
then there is a simple path from $\b$ to $\a$
\iff\  there is a long root which belongs to both $J$ and $K$.
Moreover, there is a path from $\b$ to $\a$ \iff\  there is a long root
$\s$ in $J$, and a long root $\t$ in $K$, such that $(\s|\t)\neq 0$.
\item\label{path negative} If $\a \le \b < 0$,
then there is a simple path from $\b$ to $\a$.
\end{enumerate}
\end{prop}

\proof{
We will prove (\ref{path positive}) by induction on
$m = \hauteur^\vee(\b) - \hauteur^\vee(\a)$.

\medskip

If $m = 0$, then $\b = \a$, and there is nothing to prove.

So we may assume that $m > 0$ and that the results holds
for $m - 1$. Thus $\a < \b$ and we have
$$\b - \a = \sum_{\g \in J} n_\g \g$$
where $J$ is a non-empty subset of $\D$, and the $n_\g$, $\g\in J$,
are positive integers. We have
$$(\b - \a~|~\b-\a) = \sum_{\g\in J} n_\g (\b|\g)
- \sum_{\g\in J} n_\g (\a|\g) > 0.$$
So there is a $\g$ in $J$ such that $(\b|\g) > 0$ or $(\a|\g) < 0$.
In the first case, let $\b' = s_\g(\b)$. It is a long root.
If $\g$ is long (resp. short),
then $\b' = \b - \g$ (resp. $\b' = \b - r\g$), so that
$\a \le \b' < \b$ (see Lemma \ref{caracterisation longue}).
We have $\b \elem{\g} \b'$ and $\hauteur^\vee(\b')
= \hauteur^\vee(\b) - 1$, so we can conclude by the induction
hypothesis.
The second case is similar: if $\a' = s_\g(\a) \in\Phi_\longue$,
then $\a < \a' \le \b$, $\a' \elem{\g} \a$, $\hauteur^\vee(\a')
= \hauteur^\vee(\a) + 1$, and we can conclude by the induction
hypothesis. This proves (\ref{path positive}).

Now (\ref{path negative}) follows, applying (\ref{path positive})
to $-\a$ and $-\b$ and using the symmetry $-1$.

Let us prove (\ref{path through}). If there is a long simple root $\s$
which belongs to $J$ and $K$, we have $\a \le -\s < \s \le \b$.
Using (\ref{path positive}), we find a simple path from $\b$ to $\s$,
then we have $\s \elem{\s} -\s$, and using (\ref{path negative}) we
find a simple path from $-\s$ to $\a$. So there is a simple path
from $\b$ to $\a$.

Suppose there is a long root
$\g$ in $J$, and a long root $\g'$ in $K$, such that $(\s|\t)\neq 0$.
Then either we are in the preceding case, or
there are long simple roots $\s\in J$ and $\t \in K$,
such that $\a \le -\t < \s \le \b$ and $\g = \s + \t$ is a root.
By Lemma \ref{through}, we have $\s \elem{\g} -\t$. Using
(\ref{path positive}) and (\ref{path negative}), we can find
simple paths from $\b$ to $\s$ and from $-\t$ to $\a$.
So there is a path from $\b$ to $\a$.

Now suppose there is a path from $\b$ to $\a$.
In this path, we must have a unique step of the form
$\s \elem{\g} -\t$, with $\s$ and $\t$ in $\D_\longue$.
We have $\s\in J$, $\t \in K$, and $(\s|\t) \neq 0$.
If moreover it is a simple path from $\b$ to $\a$,
then we must have $\t = -\s$. This completes the proof.
}

The preceding analysis can be used to study the length and the
reduced expressions of some elements of $W$.

\begin{prop}\label{minimal length}
Let $\a$ and $\b$ be two long roots. If $x$ is an element
of $W$ such that $x(\b) = \a$, then we have $l(x) \ge |L(\a) - L(\b)|$.

Moreover, there is an $x\in W$ such that $x(\b) = \a$ and
$l(x) = |L(\a) - L(\b)|$ \iff\ 
$\a$ and $\b$ are linked by a simple path, either
from $\b$ to $\a$, or from $\a$ to $\b$. In this case,
there is only one such $x$, and we denote it by $x_{\a\b}$.
The reduced expressions of $x_{\a\b}$ correspond bijectively
to the simple paths from $\b$ to $\a$.

If $\a \le \b \le \g$ are such that $x_{\a\b}$ and $x_{\b\g}$ are
defined, then $x_{\a\g}$ is defined, and we have $x_{\a\g} =
x_{\a\b} x_{\b\g}$ with $l(x_{\a\g}) = l(x_{\a\b}) + l(x_{\b\g})$.

The element $x_{-\a,\a}$ is defined for all $\a \in \Phi_\longue^+$,
and is equal to $s_\a$.

The element $x_{\a,\alpt}$ is defined for all $\a \in \Phi_\longue$,
and is equal to $x_\a$.
\end{prop}

\proof{
Let $(s_{\g_k},\ldots,s_{\g_1})$ be a reduced expression of $x$,
where $k = l(x)$. For $i \in \{0,~\ldots,~k\}$, let
$\b_i = s_{\g_i} \ldots s_{\g_1}(\b)$.
For each $i \in \{0,~\ldots,~k - 1\}$, we have
$|L(\b_{i + 1}) - L(\b_i)| \le 1$ by Lemma \ref{action simple}.
Then we have
$$|L(\a) - L(\b)| \le \sum_{i = 0}^{k - 1} |L(\b_{i + 1}) - L(\b_i)|
\le k = l(x)$$

If we have an equality, then all the $L(\b_{i + 1}) - L(\b_i)$ must
be of absolute value one and of the same sign, so either they are
all equal to $1$, or they are all equal to $-1$. Thus, either
we have a simple path from $\b$ to $\a$, or we have a simple path
from $\b$ to $\a$.

Suppose there is a simple path from $\b$ to $\a$.
Let $(\b_0,\ldots,\b_k)$ be a sequence of long roots,
and $(\g_1,\ldots,\g_k)$ a sequence of simple roots,
such that
$$\b = \b_0 \elem{\g_1} \b_1 \elem{\g_2} \ldots \elem{\g_k} \b_k = \a.$$
Let $x = s_{\g_k}\ldots s_{\g_1}$.
Then we have $x(\b) = \a$, and $l(x) \le k$. But we have seen that
$l(x) \ge L(\a) - L(\b) = k$. So we have equality. The case
where there is a simple path from $\a$ to $\b$ is similar.

Let $\a\in\Phi_\longue$. If $\a > 0$, we have $0 < \a \le \alpt$, so
by Proposition \ref{simple path} (\ref{path positive}), there is a simple
path from $\alpt$ to $\a$. If $\a < 0$, we have $\a < -\a \le \alpt$,
so by Proposition
\ref{simple path} (\ref{path positive}) and (\ref{path through}),
there is also a simple path from $\alpt$ to $\a$ in this case.
Let $x$ be the product of the simple reflections it involves.
Then $l(x) = L(\a)$, so $x$ is of
minimal length in $xW_{\Iti}$, and $x\in X_{\Iti}$. Thus $x = x_\a$
is uniquely determined, and $x_{\a,\alpt}$ is defined. It is equal
to $x_\a$ and is of length $L(\a)$.

Let $\a$ and $\b$ be two long roots such there is a simple path
from $\b$ to $\a$, and let $x$ be the product of the simple
reflections it involves. We have $x x_\b(\alpt) = x(\b) = \a$,
and it is of length $L(\a)$, so it is of minimal length in its coset
modulo $W_{\Iti}$. Thus $x x_\b = x_\a$, and $x = x_\a x_\b^{-1}$ is
uniquely determined. Therefore, $x_{\a\b}$ is defined and equal
to $x_\a x_\b^{-1}$. Any simple path from $\b$ to $\a$ gives
rise to a reduced expression of $x_{\a\b}$, and every reduced expression
of $x_{\a\b}$ gives rise to a simple path from $\b$ to $\a$.
These are inverse bijections.

If $\a \le \b \le \g$ are such that $x_{\a\b}$ and $x_{\b\g}$ are
defined, one can show that $x_{\a\g}$ is defined, and that
we have $x_{\a\g} = x_{\a\b} x_{\b\g}$
with $l(x_{\a\g}) = l(x_{\a\b}) + l(x_{\b\g})$,
by concatenating simple paths from $\g$ to $\b$ and from $\b$
to $\a$.

If $\a$ is a positive long root, then there is a simple path
from $\a$ to $-\a$. We can choose a symmetric path (so that the
simple reflections form a palindrome). So $x_{-\a,\a}$ is defined,
and is a reflection: it must be $s_\a$. It is of length
$L(-\a) - L(\a) = 2 \hauteur^\vee(\a) - 1$.
}

\begin{remark}\label{second proof}
\emph{We have seen in the proof that, if $\a\in\Phi_\longue^+$, then
$$l(s_\a) = l(x_{-\a,\a}) = L(-\a) - L(\a) = 2 \hauteur^\vee(\a) - 1$$
and, if $\a\in\Phi_\longue$, then
$$l(x_\a) = l(x_{\a,\alpt}) = L(\a)$$
thus we have a second proof of Lemma \ref{longueur reflexion}
and Proposition \ref{xa}.
Similarly, the formulas $x_{\a\g} = x_{\a\b} x_{\b\g}$
and $l(x_{\a\g}) = l(x_{\a\b}) + l(x_{\b\g})$, applied to
the triple $(-\a,\a,\alpt)$, give another proof of Proposition
\ref{saxa}.}
\end{remark}

\bigskip

To conclude this section, let us summarize the results which
we will use in the sequel.

\begin{theo}\label{conclu racines}
The bijection \ref{bijection} is an anti-isomorphism between the posets
$(\Phi_\longue,\preceq)$ and $(X_\Iti,\le)$
(these orders were defined at the beginning of \ref{orders}), 
and a root of level $i$ corresponds to an element of length $i$
in $X_\Iti$.

If $\b$ and $\a$ are long roots, and $\g$ is a positive root,
then we have
$$\b \elem{\g} \a \quad \textrm{ \iff\  } \quad x_\b \elem{\g} x_\a$$
(these relations have been defined at the beginning of \ref{orders}).

Moreover, in the above situation, the integer
$\partial_{\a\b} = \langle \b, \g^\vee \rangle$ is
determined as follows:
\begin{enumerate}[(i)]
\item if $\b\in\D_\longue$ and $\a\in -\D_\longue$, then
$\partial_{\a\b}$ is equal to $2$ if $\a = -\b$, and to $-1$
if $\b + (-\a)$ is a root;

\item otherwise, $\partial_{\a\b}$ is equal to $1$ if $\g$ is long,
and to $r$ if $\g$ is short (where $r = \max_{\a\in\Phi} (\a|\a)$).
\end{enumerate}
\end{theo}

If $\b$ and $\a$ are two long roots such that $L(\a) = L(\b) + 1$,
then we set $\partial_{\a\b} = 0$ if there is no simple root $\g$
such that $\b \elem{\g} \a$.

The numbers $\partial_{\a\b}$ will appear in Theorem \ref{spectral roots}
as the coefficients of the matrices of some maps appearing in the
Gysin sequence associated to the $\CM^*$-fibration $\OC_\mini \simeq
G\times_P \CM^* x_\mini$ over $G/P$, giving the cohomology of $\OC_\mini$.
By Theorem \ref{conclu racines}, these coefficients are explicitly
determined in terms of the combinatorics of the root system.

\section{Resolution of singularities, Gysin sequence}\label{res spectral}

Let us choose a maximal torus $T$ of $G$, with Lie algebra $\tG \incl \gG$.
We then denote by $X(T)$ its group of characters, and $X^\vee(T)$
its group of cocharacters. For each $\a\in\Phi$, there is a
closed subgroup $U_\a$ of $G$, and an isomorphism $u_\a : \GM_a \to U_\a$
such that, for all $t\in T$ and for all $\l\in \CM$, we have
$t u_\a(\l) t^{-1} = u_\a(\a(t) \l)$.
We are in the set-up of \ref{root systems}, with $\Phi$ equal to
the root system of $(G,T)$ in $V = X(T)\otimes_\ZM \RM$.
We denote $X(T) \times_\ZM \QM$ by $V_\QM$, and the symmetric
algebra $S(V_\QM)$ by $S$.

There is a root subspace decomposition
$$ \gG = \tG \oplus \left(\bigoplus_{\a\in\Phi} \gG_\a\right)$$
where $\gG_\a$ is the (one-dimensional) weight subspace
$\{x\in\gG~\mid~\forall t\in T,~\Ad(t).x = \a(t)x\}$.
We denote by $e_\alpha$ a non-zero vector in $\gG_\a$.
Thus we have $\gG_\a = \CM e_\a$.

Let $W = N_G(T)/T$ be the Weyl group. It acts on $X(T)$, and hence
on $V_\QM$ and $S$.

Let us now fix a Borel subgroup $B$ of $G$ containing $T$, with Lie
algebra $\bG$. This choice determines a basis $\D$, the subset of
positive roots $\Phi^+$, and the height (and dual height) function,
as in \ref{basis}, the length function $l$ as in \ref{length},
the highest root $\alpt$ and the subset $\Iti$ of $\D$ as in
\ref{highest root}, and the orders on $\Phi_\longue$ and $X_\Iti$
as in \ref{orders}.
So we can apply all the notations and results of section \ref{philg xi}.

Let $H$ be a closed subgroup of $G$, and $X$ a variety with a left
$H$-action. Then $H$ acts on $G \times X$
on the right by $(g,x).h = (gh,h^{-1}x)$.
If the canonical morphism $G\to G/H$ has local sections, then
the quotient variety $(G \times X)/H$ exists (see \cite[\S 5.5]{SPR}).
The quotient is denoted by $G\times_H X$. One has a morphism
$G\times_H X \to G/H$ with local sections, whose fibers are isomorphic
to $X$. The quotient is the \emph{fibre bundle over $G/H$ associated to
$X$}. We denote the image of $(g,x)$ in this quotient by $g *_H x$,
or simply $g * x$ if the context is clear. Note that $G$ acts on the left
on $G \times_H X$, by $g'.g *_H x = g'g *_H x$.

In \ref{G/B}, we describe the cohomology of $G/B$, both in terms Chern
classes of line bundles and in terms of fundamental classes
of Schubert varieties, and we state the Pieri formula (see
\cite{BGG, DEM, HILLER} for a description of Schubert calculus).
In \ref{parab inv}, we explain how this generalizes to the parabolic
case. In \ref{coh}, we give an algorithm to compute 
the cohomology of any line bundle minus the zero section,
on any generalized flag variety. To do this, we need the
Gysin sequence (see for example \cite{BT, HUS}, or \cite{MILNE}
in the \'etale case). In \ref{res sing}, we will see that
the computation of the cohomology of $\OC_\mini$ is a particular case.
Using the results of section \ref{philg xi}, we give a description
in terms of the combinatorics of the root system.

\subsection{Line bundles on $G/B$, cohomology of $G/B$}\label{G/B}

Let $\BC = G/B$ be the flag variety. It is a smooth projective
variety of dimension $|\Phi^+|$. The map $G \longto G/B$
has local sections (see \cite[\S 8.5]{SPR}).
If $\a$ is a character of $T$, one can lift it to $B$:
let $\CM_\a$ be the corresponding one-dimensional representation 
of $B$. We can then form the $G$-equivariant line bundle
\equat
\LC(\a) = G \times_B \CM_\a \longto G/B.
\endequat

Let $c(\a)\in H^2(G/B,\ZM)$ denote the first Chern class of $\LC(\a)$.
Then $c : X(T) \longto H^2(G/B,\ZM)$ is a morphism of $\ZM$-modules.
It extends to a morphism of $\QM$-algebras, which we still
denote by $c : S \longto H^*(G/B,\QM)$. The latter is surjective
and has kernel $\IC$, where $\IC$ is the ideal of $S$ generated
by the $W$-invariant homogeneous elements in $S$ of positive degree.
So it induces an isomorphism of $\QM$-algebras
\equat
\ov{c}:S/\IC \simeq H^*(G/B,\QM)
\endequat
which doubles degrees.

The algebra $S/\IC$ is called the coinvariant algebra.
As a representation of $W$, it is isomorphic to the regular
representation. We also have an action of $W$ on $H^*(G/B,\QM)$,
because $G/B$ is homotopic to $G/T$, and $W$ acts on the right on $G/T$ by
the formula $gT.w = gnT$, where $n\in N_G(T)$ is a 
representative of $w \in W$, and $g \in G$. One can show
that $\ov c$ commutes with the actions of $W$.

On the other hand, we have the Bruhat decomposition \cite[\S 8.5]{SPR}
\equat
G/B = \bigsqcup_{w\in W} C(w)
\endequat
where the $C(w) = BwB/B \simeq \CM^{l(w)}$ are the Schubert cells.
Their closures are the Schubert varieties $S(w) = \ov{C(w)}$.
Thus the cohomology of $G/B$ is concentrated in even degrees,
and $H^{2i}(G/B,\ZM)$ is free with basis $(Y_w)_{w\in W^i}$,
where $Y_w$ is the cohomology class of the Schubert variety $S(w_0 w)$
(which is of codimension $l(w)=i$).
The object of Schubert calculus is to describe
the multiplicative structure of $H^*(G/B,\ZM)$ in these
terms (see \cite{BGG, DEM, HILLER}).
We will only need the following result (known
as the Pieri formula, or Chevalley formula):
if $w\in W$ and $\a\in X(T)$, then 
\equat\label{pieri}
c(\a)\ .\  Y_w\quad = \sum_{w\elem{\g}w'}
                   \left<w(\a), \g^\vee\right>  Y_{w'}
\endequat

\subsection{Parabolic invariants}\label{parab inv}

Let $I$ be a subset of $\D$.
Let $P_I$ be the parabolic subgroup of $G$ containing $B$ corresponding
to $I$. It is generated by $B$ and the subgroups $U_{-\a}$, for
$\a\in I$. Its unipotent radical $U_{P_I}$ is generated by the $U_\a$,
$\a\in \Phi^+\setminus\Phi_I^+$. And it has a Levi complement
$L_I$, which is generated by $T$ and the $U_\a$, $\a\in\Phi_I$.
One can generalize the preceding constructions to the parabolic case.

If $\a\in X(T)^{W_I}$ (that is, if $\a$ is a character orthogonal
to $I$), then we can form the $G$-equivariant line bundle
\equat
\LC_I(\a) = G \times_{P_I} \CM_\alpha \longto G/{P_I}
\endequat
because the character $\a$ of $T$, invariant by $W_I$, can be extended to $L_I$
and lifted to ${P_I}$.

We have a surjective morphism $q_I : G/B \longto G/P_I$, which induces
an injection
\[
q_I^* : H^*(G/P_I,\ZM) \hookrightarrow H^*(G/B,\ZM)
\]
in cohomology,
which identifies $H^*(G/P_I,\ZM)$ with $H^*(G/B,\ZM)^{W_I}$.

The isomorphism $\ov{c}$ restricts to
\equat
(S/\IC)^{W_I} \simeq H^*(G/{P_I},\QM)
\endequat

We have cartesian square
\[
\xymatrix{
\LC(\a) \ar[r] \ar[d] & \LC_I(\a) \ar[d]\\
G/B \ar[r]_{q_I} & G/P_I
}
\]

That is, the pullback by $q_I$  of $\LC_I(\a)$ is $\LC(\a)$.
By functoriality of Chern classes, we have $q_I^*(c_I(\a)) = c(\a)$.

We still have a Bruhat decomposition
\equat
G/P_I = \bigsqcup_{w\in X_I} C_I(w)
\endequat
where $C_I(w) =  BwP_I/P_I \simeq \CM^{l(w)}$ for $w$ in $X_I$.
We note
\[
\text{
$S_I(w) = \ov{C_I(w)}$ and
$Y_{I,w} = [\ov{Bw_0 wP_I/P_I}] = [\ov{Bw_0 w w_I P_I/P_I}]
= [S_I(w_0 w w_I)]$
for $w$ in $X_I$}
\]
Note that
\equat
\text{if $w$ is in $X_I$, then $w_0 w w_I$ is also in $X_I$}
\endequat
since for any root $\b$ in $\Phi_I^+$, we have
$w_I(\b) \in \Phi_I^-$, hence $ww_I(\b)$ is also negative, and thus
$w_0 w w_I(\b)$ is positive. Moreover, we have
\equat
\text{if $w \in X_I$, then
$
l(w_0 w w_I)
= l(w_0) - l(w w_I)
= l(w_0)-l(w_I)-l(w)
= \dim G/P_I - l(w)$}
\endequat
so that $Y_{I,w} \in H^{2l(w)}(G/P_I,\ZM)$.
We have $q_I^*(Y_{I,w}) = Y_w$.

The cohomology of $G/P_I$ is concentrated in even degrees,
and $H^{2i}(G/P_I,\ZM)$ is free with basis $(Y_{I,w})_{w\in X_I^i}$.
The cohomology ring $H^*(G/P_I,\ZM)$ is identified \emph{via} $q_I^*$ to
a subring of $H^*(G/B,\ZM)$,
so the Pieri formula can now be written as follows.
If $w\in X_I$ and $\a\in X(T)^{W_I}$, then we have
\equat\label{pieri p}
c_I(\a)\ .\  Y_{I,w}
\quad
=
\sum_{w\elem{\g}w'\in X_I}
  \left<w(\a), \g^\vee\right>  Y_{I,w'}
\endequat

\subsection{Cohomology of a $\CM^*$-fiber bundle on $G/P_I$}\label{coh}

Let $I$ be a subset of $\D$, and $\a$ be a $W_I$-invariant character
of $T$. Let us consider
\equat\label{fib}
\LC^*_I(\a) = G \times_{P_I} \CM^*_\alpha \longto G/{P_I},
\endequat
that is, the line bundle $\LC_I(\a)$ minus the zero section.
In the sequel, we will have to calculate the cohomology
of $\LC^*_\Iti(\alpt)$, but we can explain how to calculate
the cohomology of $\LC^*_I(\a)$ for any given $I$ and $\a$
(the point is that the answer for the middle cohomology will turn out to be
nicer in our particular case, thanks to the results of section
\ref{philg xi}).

We have the Gysin exact sequence
\[
H^{n-2}(G/P_I,\ZM) \elem{c_I(\a)} H^n(G/P_I,\ZM)
\longto H^n(\LC^*_I(\a),\ZM) \longto 
H^{n-1}(G/P_I,\ZM) \elem{c_I(\a)} H^{n+1}(G/P_I,\ZM)
\]
where $c_I(\a)$ means the multiplication by $c_I(\a)$,
so we have a short exact sequence
\[
0\longto \Coker\ (c_I(\a):H^{n-2}\to H^n)
\longto H^n(\LC^*_I(\a),\ZM)
\longto \Ker\ (c_I(\a):H^{n-1}\to H^{n+1})
\longto 0
\]
where $H^j$ stands for $H^j(G/P_I,\ZM)$.

Moreover, by the hard Lefschetz theorem, $c_I(\a):\QM\otimes_\ZM H^{n-2}
\to \QM\otimes_\ZM H^n$
is injective for $n \leqslant d_I = \dim \LC^*_I(\a) = \dim G/P_I + 1$,
and surjective for $n \geqslant d_I$. By the way, we see that we could
immediately determine the rational cohomology of $\OC_\mini$,
using only the results in this paragraph and the cohomology of
$G/P_I$.

But we can say more.
The cohomology of $G/P_I$ is free and concentrated in even degrees.
In fact, $c_I(\a): H^{n-2} \to H^n$ is injective for $n \leqslant d_I$,
and has free kernel and finite cokernel for $n \geqslant d_I$.

We have 
\equat
\text{if $n$ is even, then } H^n(\LC^*_I(\a),\ZM)\simeq
\Coker\ (c_I(\a):H^{n-2}\to H^n)
\endequat
which is finite for $n \geqslant d_I$, and
\equat
\text{if $n$ is odd, then } H^n(\LC^*_I(\a),\ZM)\simeq
\Ker\ (c_I(\a):H^{n-1}\to H^{n+1})
\endequat
which is free (it is zero if $n \leqslant d_I - 1$).

Thus all the cohomology of $\LC^*_I(\a)$ can be
explicitly computed, thanks to the results of \ref{parab inv}.

\subsection{Resolution of singularities}\label{res sing}

Let $\Iti$ be the subset of $\D$ defined in \ref{Iti}.
There is a resolution of singularities (see for example the introduction
of \cite{KP2})

\equat
\fonctio{
\LC_\Iti(\alpt) = G \times_{P_\Iti} \CM_\alpt }{
\overline{\OC_\mini} = \OC_\mini \cup \{0\} }{
g*\l }{ \Ad g.(\l . e_\alpt)}
\endequat

It is the one mentioned in the introduction, with $P = P_\Iti$ and
$x_\mini = e_\alpt$.
It induces an isomorphism

\equat
\LC^*_\Iti(\alpt) = G \times_{P_\Iti} \CM^*_\alpt \simeq \OC_\mini
\endequat

Set $d = d_\Iti = 2 h^\vee -2$. For all integers $j$, let
$H^j$ denote $H^j(G/P_\Iti,\ZM)$.
For $\a\in\Phi_\longue$, we have $x_\a \in X_\Iti$, so
$Z_\a := Y_{\Iti,x_\a}$  is an element of $H^{2i}(G/P_\Iti,\ZM)$,
where $i = l(x_\a) = L(\a)$.
Then $H^*(G/P_\Iti,\ZM)$ is concentrated in even degrees, and
$H^{2i}(G/P_\Iti,\ZM)$ is free with basis $(Z_\a)_{\a\in \Phi_\longue^i}$.
Combining Theorem \ref{conclu racines} and the analysis of \ref{coh},
we get the following description of the cohomology of $\OC_\mini$.

\begin{theo}\label{spectral roots}
We have
\equat
\text{if $n$ is even, then } H^n(\OC_\mini,\ZM)\simeq
\Coker\ (c_\Iti(\alpt):H^{n-2}\to H^n)
\endequat
which is finite for $n \geqslant d$, and
\equat
\text{if $n$ is odd, then } H^n(\LC^*_I(\a),\ZM)\simeq
\Ker\ (c_\Iti(\alpt):H^{n-1}\to H^{n+1})
\endequat
which is free (it is zero if $n \leqslant d - 1$).

Moreover, if $\b\in\Phi_\longue^i$, then we have
\equat
c_\Iti(\alpt).Z_\b = 
\sum_{\b\elem{\g}\a} \left<\b, \g^\vee\right>  Z_\a
= \sum_{\a\in\Phi_\longue^{i + 1}} \partial_{\a\b} Z_\b
\endequat
where the $\partial_{\a\b}$ are the integers defined in
Theorem \ref{conclu racines}.
\end{theo}

As a consequence, we obtain the following results.

\begin{theo}\label{middle}
\begin{enumerate}[(i)]
\item
The middle cohomology of $\OC_\mini$ is given by
$$H^{2h^\vee - 2}(\OC_\mini, \ZM) \simeq P^\vee({\Phi'})/Q^\vee({\Phi'})$$
where $\Phi'$ is the sub-root system of $\Phi$ generated by $\D_\longue$,
and $P^\vee({\Phi'})$ (resp. $Q^\vee({\Phi'})$) is its coweight lattice
(resp. its coroot lattice).

\item
The rest of the cohomology of $\OC_\mini$ is as described
in section \ref{case by case}. In particular,
if $\ell$ is a good prime for $G$, then there is no $\ell$-torsion
in the rest of the cohomology of $\OC_\mini$.
\end{enumerate}
\end{theo}

\proof{
The map
$c_\Iti(\alpt) : H^{2 h^\vee - 4} \longto H^{2 h^\vee -2}$
is described as follows. By Theorem \ref{spectral roots},
the cohomology group $H^{2 h^\vee - 4}$ is free with basis
$(Z_\b)_{\b\in \D_\longue}$ (the long roots
of level $h^\vee - 2$ are of the long roots dual height $1$, so they are
the long simple roots). Similarly,
$H^{2 h^\vee - 4}$ is free with basis
$(Z_{-\a})_{\a\in \D_\longue}$. Besides, the matrix of
$c_\Iti(\alpt):H^{2 h^\vee - 4} \longto H^{2 h^\vee -2}$
in these bases is $(\partial_{-\a,\b})_{\a,\b\in\D_\longue}$.
We have $\partial_{-\a,\a} = 2$ for $\a\in\D_\longue$, and
for distinct $\a$ and $\b$ we have $\partial_{-\a,\b} = 1$ if
$\a + \b$ is a (long) root, $0$ otherwise.

Thus the matrix of $c_\Iti(\alpt):H^{2 h^\vee - 4} \longto H^{2 h^\vee -2}$
is the Cartan matrix of $\Phi'$ without minus signs.
This matrix is equivalent to the Cartan matrix of $\Phi'$:
since the Dynkin diagram of $\Phi'$ is a tree, it is bipartite.
We can write $\D_\longue = J \cup K$, where no element of $J$
is linked to an element of $K$ in the Dynkin diagram of $\Phi'$.
If we replace the $Z_{\pm \a}$, $\a \in J$, by their opposites, 
then the matrix of $c_\Iti(\alpt) : H^{2 h^\vee - 4} \longto H^{2 h^\vee -2}$
becomes the Cartan matrix of $\Phi'$.

Now, the Cartan matrix of $\Phi'$ is transposed to the matrix
of the inclusion of $Q(\Phi')$ in $P(\Phi')$ in the bases
$\D_\longue$ and $(\varpi_\a)_{\a\in\D_\longue}$ (see
\cite[Chap. VI, \S 1.10]{BOUR456}), so it is in fact
the matrix of the inclusion of $Q^\vee(\Phi')$ in $P^\vee(\Phi')$
in the bases $(\b^\vee)_{\b\in\D_\longue}$ and
$(\varpi_{\a^\vee})_{\a\in\D_\longue}$.

The middle cohomology group $H^{2h^\vee - 2}(\OC_\mini, \ZM)$
is isomorphic to the cokernel of the map
$c_\Iti(\alpt):H^{2 h^\vee - 4} \longto H^{2 h^\vee -2}$.
This proves $(i)$.

Part $(ii)$ follows from a case-by-case analysis which will be done
in section \ref{case by case}.
}

\begin{remark}\label{poincare}
\emph{Besides, we have $\partial_{\a\b} = \partial_{-\b,-\a}$, so the
maps ``multiplication by $c_\Iti(\alpt)$''
in complementary degrees are transposed to each other.
This accounts for the fact that $\OC_\mini$ satisfies Poincar\'e duality,
since $\OC_\mini$ is homeomorphic to $\RM^+_*$ times a smooth compact
manifold of (real) dimension $2 h^\vee - 5$ (since we deal with
integral coefficients, one should take the derived dual for the Poincar\'e
duality).}
\end{remark}

\begin{remark}\label{rat}
\emph{
For the first half of the rational cohomology of $\OC_\mini$,
we find
\[
\bigoplus_{i=1}^k \QM[-2(d_i - 2)]
\]
where $k$ is the number of long simple roots, and
$d_1 \leqslant \ldots \leqslant d_k \leqslant \ldots \leqslant d_n$
are the degrees of $W$ ($n$ being the total number of simple roots).
This can be observed case by case, or related to the corresponding
Springer representation. The other half is determined by Poincar\'e
duality.
}
\end{remark}

\section{Case-by-case analysis}\label{case by case}

In the preceding section, we have explained how to compute the cohomology
of the minimal class in any given type in terms of root systems,
and we found a description of
the middle cohomology with a general proof.
However, for the rest of the cohomology, we need a case-by-case analysis.
It will appear that the primes dividing the torsion of the rest of the
cohomology are bad.
We have no \emph{a priori} explanation for this fact. Note that,
for the type $A$, we have an alternative method, which will be explained
in the next section.

For all types, first we give the Dynkin diagram, to fix the numbering
$(\a_i)_{1 \le i \le r}$ of the vertices, where $r$ denotes
the semisimple rank of $\gG$,
and to show the part $\Iti$ of $\D$ (see \ref{Iti}).
The corresponding vertices are represented in black. They are
exactly those that are not linked to the additional vertex in the
extended Dynkin diagram.

Then we give a diagram whose vertices are the positive long roots;
whenever $\b \elem{\g} \a$, we put an edge between $\b$ (above) and
$\a$ (below), and the multiplicity of the edge is equal to
$\partial_{\a\b} = \langle \b, \g^\vee \rangle$. In this diagram, the long root
$\sum_{i = 1}^r n_i \a_i$ (where the $n_i$ are non-negative integers)
is denoted by $n_1\ldots n_r$.
The roots in a given line appear in lexicographic order.

For $1 \leqslant i \leqslant d - 1$, let $\DC_i$ be the
matrix of the map $c_\Iti(\alpt):H^{2i - 2} \to H^{2i}$ 
in the bases $\Phi_\longue^{i - 1}$ and $\Phi_\longue^{i}$ (the roots
being ordered in lexicographic order, as in the diagram).
We give the matrices $\DC_i$
for $i = 1 \ldots h^\vee - 2$.
The matrix $\DC_{h^\vee - 1}$ is equal to
the Cartan matrix without minus signs of the root system $\Phi'$
(corresponding to $\D_\longue$).
The last matrices can be deduced
from the first ones by symmetry, since (by Remark \ref{poincare}) we have
$\DC_{d - i} = {}^t \DC_i$.

Then we give the cohomology of the minimal class with $\ZM$ coefficients
(one just has to compute the elementary divisors of the matrices $\DC_i$).

It will be useful to introduce some notation for the matrices in classical
types. Let $k$ be an integer. We set
\[
M(k) = 
\begin{pmatrix}
1 & 0      & \ldots & 0 & 0\\
1 & 1      & \ddots & \vdots & \vdots\\
0 & 1      & \ddots & 0 & 0\\
\vdots & \ddots & \ddots & 1 & 0\\
0 & \ldots & 0      & 1 & 1
\end{pmatrix}
\hspace{3cm}
N(k) = 
\begin{pmatrix}
1 & 0      & \ldots & 0\\
1 & 1      & \ddots & \vdots\\
0 & 1      & \ddots & 0\\
\vdots & \ddots & \ddots & 1\\
0 & \ldots & 0      & 1
\end{pmatrix}
\]
where $M(k)$ is a square matrix of size $k$, and $N(k)$ is
of size $(k + 1) \times k$.

Now let $k$ and $l$ be non-negative integers. For $i$ and $j$
any integers, we define a $k \times l$ matrix $E_{i,j}(k,l)$ as follows.
If $(i, j)$ is not in the range $[1,k] \times [1,l]$, then
we set $E_{i,j}(k,l) = 0$, otherwise it will denote the
$k\times l$ matrix whose only non-zero entry is a $1$ in the intersection
of line $i$ and column $j$. If the size of the matrix is clear from the
context, we will simply write $E_{i,j}$.

First, the calculations of the elementary divisors of the matrices 
$\DC_i$ were
done with GAP3 (see \cite{GAP}). We used the data on roots systems 
of the package CHEVIE. But actually, all the calculations can be
done by hand.

\subsection{Type $A_{n - 1}$}

\begin{center}
\begin{picture}(280,30)(-35,0)
\put(  5, 10){\circle{10}}
\put( 45, 10){\circle*{10}}
\put( 85, 10){\circle*{10}}
\put(165, 10){\circle*{10}}
\put(205, 10){\circle*{10}}
\put(245, 10){\circle{10}}
\put( 10, 10){\line(1,0){30}}
\put( 50, 10){\line(1,0){30}}
\put( 90, 10){\line(1,0){20}}
\put(120,9.5){$\ldots$}
\put(140, 10){\line(1,0){20}}
\put(170, 10){\line(1,0){30}}
\put(210, 10){\line(1,0){30}}
\put(  0, 20){$\a_1$}
\put( 40, 20){$\a_2$}
\put( 80, 20){$\a_3$}
\put(155, 20){$\a_{n - 3}$}
\put(195, 20){$\a_{n - 2}$}
\put(235, 20){$\a_{n - 1}$}
\end{picture}
\end{center}

We have $h = h^\vee = n$ and $d = 2n - 2$.

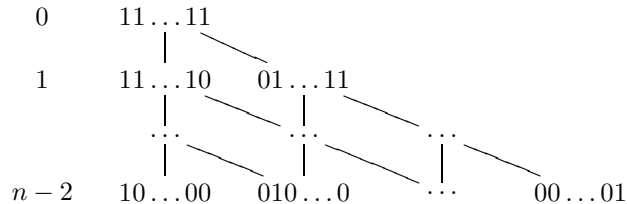
\[
\xymatrix @=.4cm{
0 &
11\ldots 11 \ar@{-}[d] \ar@{-}[dr]\\
1 &
11\ldots 10 \ar@{-}[d] \ar@{-}[dr]&
01\ldots 11 \ar@{-}[d] \ar@{-}[dr]&
{}\phantom{01\ldots 10}
\\
&
\ldots \ar@{-}[d] \ar@{-}[dr]&
\ldots \ar@{-}[d] \ar@{-}[dr]&
\ldots \ar@{-}[d] \ar@{-}[dr]\\
n - 2 &
10\ldots 00 &
010\ldots 0 &
\ldots &
00\ldots 01
}
\]

The odd cohomology of $G/P_\Iti$ is zero, and we have
\[
H^{2i}(G/P_\Iti) =
\begin{cases}
\ZM^{i + 1}
& \text{if $0 \le i \le n - 2$}\\
\ZM^{2n - 2 - i}
& \text{if $n - 1 \le i \le 2n - 3$}\\
0 & \text{otherwise}
\end{cases}
\]

For $1 \le i \le n - 2$, we have $\DC_i = N(i)$; the cokernel is isomorphic
to $\ZM$. We have
\[
\DC_{n - 1} =
\begin{pmatrix}
2      & 1      & 0      & \ldots & 0      \\
1      & 2      & 1      & \ddots & \vdots \\
0      & \ddots & \ddots & \ddots & 0      \\
\vdots & \ddots & 1      & 2      & 1      \\
0      & \ldots & 0      & 1      & 2
  \end{pmatrix}
\]
Its cokernel is isomorphic to $\ZM/n$. The last matrices are
transposed to the first ones, so the corresponding
maps are surjective. From this, we deduce the cohomology
of $\OC_\mini$ in type $A_{n - 1}$. We will see another method in
section \ref{a bis}.

\[
H^i(\OC_\mini, \ZM) =
\begin{cases}
\ZM & \text{if $0 \le i \le 2n -4$ and $i$ is even,}\\
& \text{or $2n - 1 \le i \le 4n - 5$ and $i$ is odd}\\
\ZM/n & \text{if $i = 2n - 2$}\\
0 & \text{otherwise}
\end{cases}
\]

\subsection{Type $B_n$}

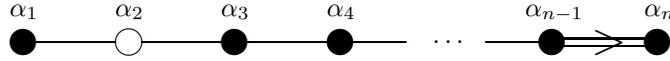
\begin{center}
\begin{picture}(280,30)(-35,0)
\put(  5, 10){\circle*{10}}
\put( 45, 10){\circle{10}}
\put( 85, 10){\circle*{10}}
\put(125, 10){\circle*{10}}
\put(205, 10){\circle*{10}}
\put(245, 10){\circle*{10}}
\put( 10, 10){\line(1,0){30}}
\put( 50, 10){\line(1,0){30}}
\put( 90, 10){\line(1,0){30}}
\put(130, 10){\line(1,0){20}}
\put(160,9.5){$\ldots$}
\put(180, 10){\line(1,0){20}}
\put(209, 11.5){\line(1,0){32}}
\put(209,  8.5){\line(1,0){32}}
\put(220,  5.5){\LARGE{$>$}}
\put(  0, 20){$\a_1$}
\put( 40, 20){$\a_2$}
\put( 80, 20){$\a_3$}
\put(120, 20){$\a_4$}
\put(195, 20){$\a_{n - 1}$}
\put(240, 20){$\a_n$}
\end{picture}
\end{center}

We have $h = 2n$, $h^\vee = 2n - 1$, and $d = 4n - 4$.

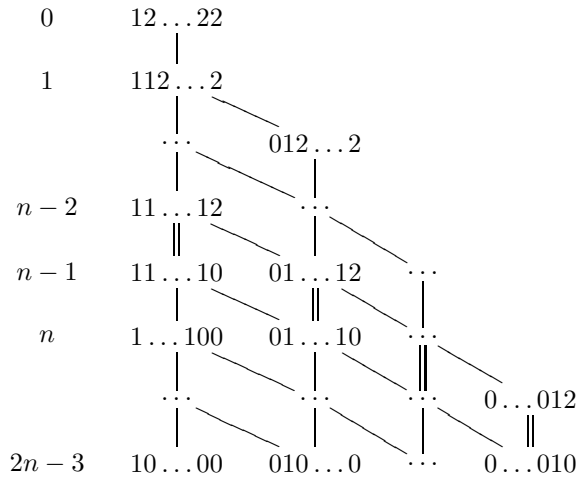
\[
\xymatrix @=.4cm{
0 & 
12\ldots 22 \ar@{-}[d]\\
1 & 
112\ldots 2 \ar@{-}[d] \ar@{-}[dr]\\
& 
\ldots \ar@{-}[d] \ar@{-}[dr]&
012\ldots 2 \ar@{-}[d]\\
n - 2 & 
11\ldots 12 \ar@{=}[d] \ar@{-}[dr]&
\ldots \ar@{-}[d] \ar@{-}[dr]&\\
n - 1 &
11\ldots 10 \ar@{-}[d] \ar@{-}[dr]&
01\ldots 12 \ar@{=}[d] \ar@{-}[dr]&
\ldots \ar@{-}[d]\\
n &
1\ldots 100 \ar@{-}[d] \ar@{-}[dr]&
01\ldots 10 \ar@{-}[d] \ar@{-}[dr]&
\ldots \ar@{=}[d] \ar@{-}[dr]\\
&
\ldots \ar@{-}[d] \ar@{-}[dr]&
\ldots \ar@{-}[d] \ar@{-}[dr]&
\ldots \ar@{-}[d] \ar@{-}[dr]&
0\ldots 012 \ar@{=}[d]\\
2n - 3 &
10\ldots 00&
010\ldots 0&
\ldots&
0\ldots 010
}
\]

There is a gap at each even line (the length of the line
increases by one). The diagram can be a little bit misleading if
$n$ is even: in that case, there is a gap at the line $n - 2$.
Let us now describe the matrices $\DC_i$.

First suppose $1 \le i \le n - 2$. If $i$ is odd, then we have
$\DC_i = M\left(\tfrac{i+1}{2}\right)$ (an isomorphism); 
if $i$ is even, then we have
$\DC_i = N\left(\tfrac{i}{2}\right)$ and the cokernel is isomorphic to $\ZM$.

Now suppose $n - 1 \le i \le 2n - 3$. If $i$ is odd,
then we have
$\DC_i = M\left(\tfrac{i+1}{2}\right) +  E_{i + 2 - n, i + 2 - n}$
and the cokernel is isomorphic to $\ZM/2$.
If $i$ is even, then we have
$\DC_i = N\left(\tfrac{i}{2}\right) + E_{i + 2 - n, i + 2 - n}$
and the cokernel is isomorphic to $\ZM$.

The long simple roots generate a root system of type $A_{n - 1}$.
Thus the matrix $\DC_{2n - 2}$ is the Cartan matrix without minus
signs of type $A_{n - 1}$, which has cokernel $\ZM/n$.

So the cohomology of $\OC_\mini$ is described as follows.

\[ H^i(\OC_\mini, \ZM) \simeq
\begin{cases}
\ZM & \text{if $0 \le i \le 4n - 8$ and $i \equiv 0 \mod 4$,}\\
& \text{or $4n - 1 \le i \le 8n  - 9$ and $i \equiv -1 \mod 4$}\\
\ZM/2 & \text{$2n - 2 \le i \le 6n - 6$ and $i \equiv 2 \mod 4$}\\
\ZM/n & \text{if $i = 4n - 4$}\\
0 & \text{otherwise}
\end{cases}
\]

\subsection{Type $C_n$}

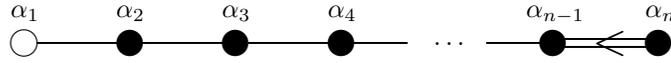
\begin{center}
\begin{picture}(280,30)(-35,0)
\put(  5, 10){\circle{10}}
\put( 45, 10){\circle*{10}}
\put( 85, 10){\circle*{10}}
\put(125, 10){\circle*{10}}
\put(205, 10){\circle*{10}}
\put(245, 10){\circle*{10}}
\put( 10, 10){\line(1,0){30}}
\put( 50, 10){\line(1,0){30}}
\put( 90, 10){\line(1,0){30}}
\put(130, 10){\line(1,0){20}}
\put(160,9.5){$\ldots$}
\put(180, 10){\line(1,0){20}}
\put(209, 11.5){\line(1,0){32}}
\put(209,  8.5){\line(1,0){32}}
\put(220,  5.5){\LARGE{$<$}}
\put(  0, 20){$\a_1$}
\put( 40, 20){$\a_2$}
\put( 80, 20){$\a_3$}
\put(120, 20){$\a_4$}
\put(195, 20){$\a_{n - 1}$}
\put(240, 20){$\a_n$}
\end{picture}
\end{center}

We have $h = 2n$, $h^\vee = n + 1$, and $d = 2n$.
The root system $\Phi'$ is of type $A_1$. Its Cartan matrix is $(2)$.

\[
\xymatrix @=.5cm{
0 &
22\ldots 21 \ar@{=}[d]\\
1 &
02\ldots 21 \ar@{=}[d]\\
&
\vdots      \ar@{=}[d]\\
n - 2 &
0\ldots 021 \ar@{=}[d]\\
n - 1 &
00\ldots 01
}
\]

The matrices $\DC_i$ are all equal to $(2)$.

\[
H^i(\OC_\mini, \ZM) \simeq
\begin{cases}
\ZM & \text{if $i = 0$ or $4n - 1$}\\
\ZM/2 & \text{if $2 \le i \le 4n - 2$ and $i$ is even}\\
0 & \text{otherwise}
\end{cases}
\]

\subsection{Type $D_n$}

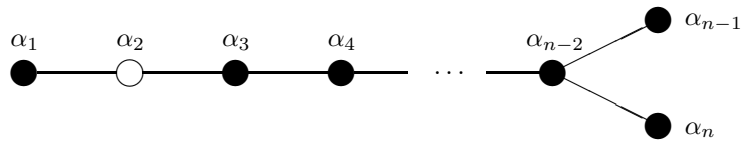
\begin{center}
\begin{picture}(280,80)(-35,-20)
\put(  5, 10){\circle*{10}}
\put( 45, 10){\circle{10}}
\put( 85, 10){\circle*{10}}
\put(125, 10){\circle*{10}}
\put(205, 10){\circle*{10}}
\put(245, 30){\circle*{10}}
\put(245,-10){\circle*{10}}
\put( 10, 10){\line(1,0){30}}
\put( 50, 10){\line(1,0){30}}
\put( 90, 10){\line(1,0){30}}
\put(130, 10){\line(1,0){20}}
\put(160,9.5){$\ldots$}
\put(180, 10){\line(1,0){20}}
\put(209, 12){\line(2,1){32}}
\put(209,  8){\line(2,-1){32}}
\put(  0, 20){$\a_1$}
\put( 40, 20){$\a_2$}
\put( 80, 20){$\a_3$}
\put(120, 20){$\a_4$}
\put(195, 20){$\a_{n - 2}$}
\put(255, 28){$\a_{n - 1}$}
\put(255,-13){$\a_n$}
\end{picture}
\end{center}

We have $h = h^\vee = 2n - 2$, and $d = 4n -6$.
We have
\[
\DC_{2n-3} =
\begin{pmatrix}
2      & 1      & 0      & \ldots & 0      & 0\\
1      & 2      & \ddots & \ddots & \vdots & \vdots\\
0      & \ddots & \ddots & 1      & 0      & 0\\
\vdots & \ddots & 1      & 2      & 1      & 1\\
0      & \ldots & 0      & 1      & 2      & 0\\
0      & \ldots & 0      & 1      & 0      & 2
\end{pmatrix}
\]
Its cokernel is $(\ZM/2)^2$ when $n$ is even, $\ZM/4$ when $n$ is odd.

As in the $B_n$ case, the reader should be warned that there is a gap
at line $n - 4$ if $n$ is even. Besides, not all dots are meaningful.
The entries $0\ldots 01211$ and $00\ldots 0111$ are on the right
diagonal, but usually they are not on the lines $n - 1$ and $n$.

{\scriptsize
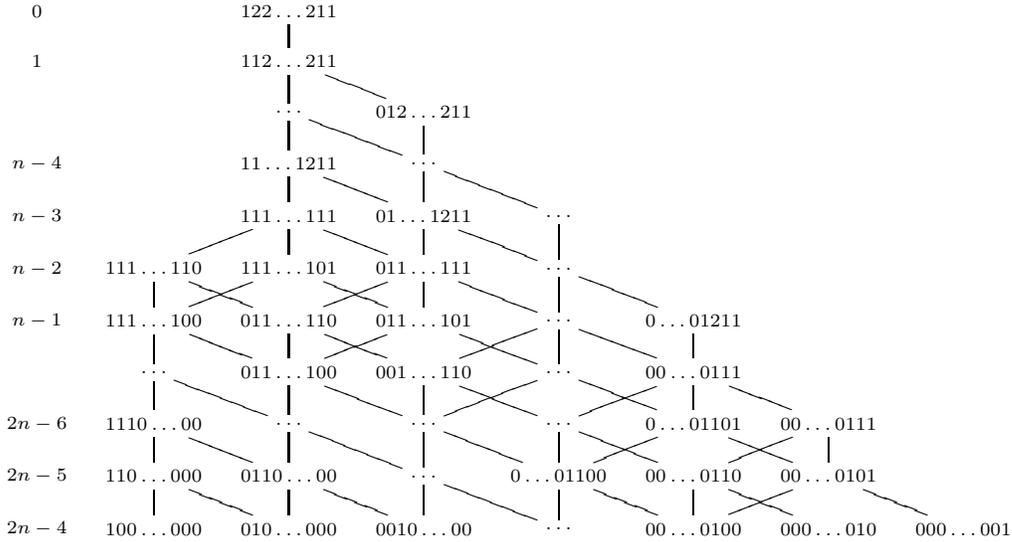
\[
\xymatrix @=.3cm{
0&
&122\ldots 211 \ar@{-}[d]\\
1&
&112\ldots 211 \ar@{-}[d] \ar@{-}[dr]\\
&
&\ldots \ar@{-}[d] \ar@{-}[dr]&
012\ldots 211  \ar@{-}[d]\\
n - 4&
&11\ldots 1211 \ar@{-}[d] \ar@{-}[dr]&
\ldots \ar@{-}[d] \ar@{-}[dr]\\
n - 3&
&111\ldots 111 \ar@{-}[dl] \ar@{-}[d] \ar@{-}[dr]&
01\ldots 1211 \ar@{-}[d] \ar@{-}[dr]&
\ldots \ar@{-}[d]\\
n - 2&
111\ldots 110 \ar@{-}[d] \ar@{-}[dr]&
111\ldots 101 \ar@{-}[dl] \ar@{-}[dr]&
011\ldots 111 \ar@{-}[dl] \ar@{-}[d] \ar@{-}[dr]&
\ldots \ar@{-}[d] \ar@{-}[dr]\\
n - 1&
111\ldots 100 \ar@{-}[d] \ar@{-}[dr]&
011\ldots 110 \ar@{-}[d] \ar@{-}[dr]&
011\ldots 101 \ar@{-}[dl] \ar@{-}[dr]&
\dots \ar@{-}[dl] \ar@{-}[d] \ar@{-}[dr]&
0\ldots 01211  \ar@{-}[d]\\
&
\ldots \ar@{-}[d] \ar@{-}[dr]&
011\ldots 100 \ar@{-}[d] \ar@{-}[dr]&
001\ldots 110 \ar@{-}[d] \ar@{-}[dr]&
\ldots \ar@{-}[dl] \ar@{-}[dr]&
00\ldots 0111 \ar@{-}[dl] \ar@{-}[d] \ar@{-}[dr]\\
2n - 6&
1110\ldots 00 \ar@{-}[d] \ar@{-}[dr]&
\ldots \ar@{-}[d] \ar@{-}[dr]&
\ldots \ar@{-}[d] \ar@{-}[dr]&
\ldots \ar@{-}[d] \ar@{-}[dr]&
0\ldots 01101 \ar@{-}[dl] \ar@{-}[dr]&
00\ldots 0111 \ar@{-}[dl] \ar@{-}[d]\\
2n - 5&
110\ldots 000 \ar@{-}[d] \ar@{-}[dr]&
0110\ldots 00 \ar@{-}[d] \ar@{-}[dr]&
\ldots \ar@{-}[d] \ar@{-}[dr]&
0\ldots 01100 \ar@{-}[d] \ar@{-}[dr]&
00\ldots 0110 \ar@{-}[d] \ar@{-}[dr]&
00\ldots 0101 \ar@{-}[dl] \ar@{-}[dr]\\
2n - 4&
100\ldots 000&
010\ldots 000&
0010\ldots 00&
\ldots&
00\ldots 0100&
000\ldots 010&
000\ldots 001
}
\]
}

First suppose $i \le i \le n - 3$. We have
\[
\DC_i =
\begin{cases}
M\left(\tfrac{i + 1}{2}\right)
& \text{if $i$ is odd}\\
N\left(\tfrac{i}{2}\right)
& \text{if $i$ is even}
\end{cases}
\]
Then the cokernel is zero if $i$ is odd, $\ZM$ if $i$ is even.

Let $V$ be the $1 \times \tfrac{n - 1}{2}$ matrix $(1,0,\ldots,0)$.
We have
\[\DC_{n - 2} =
\begin{cases}
\begin{pmatrix}
V\\
N\left(\frac{n - 2}{2}\right)
\end{pmatrix}
& \text{if $n$ is even}\\
\\
\begin{pmatrix}
V\\
M\left(\frac{n - 1}{2}\right)
\end{pmatrix}
& \text{if $n$ is odd}\\
\end{cases}
\]
The cokernel is $\ZM^2$ if $n$ is even, $\ZM$ if $i$ is odd.

Now suppose $n - 1 \le i \le 2n - 4$. We have
\[
\DC_i =
\begin{cases}
M\left(\tfrac{i + 3}{2}\right)
+ E_{i + 2 - n, i + 3 - n}
- E_{i + 3 - n, i + 3 - n}
+ E_{i + 3 - n, i + 4 - n}
& \text{if $i$ is odd}\\
N\left(\tfrac{i + 2}{2}\right)
+ E_{i + 2 - n, i + 3 - n}
- E_{i + 3 - n, i + 3 - n}
+ E_{i + 3 - n, i + 4 - n}
& \text{if $i$ is even}\\
\end{cases}
\]
Then the cokernel is $\ZM/2$ if $i$ is odd, $\ZM$ if $i$ is even.

\[
H^i(\OC_\mini, \ZM) \simeq
\begin{cases}
\ZM & \text{if $0 \le i \le 4n - 8$ and $i \equiv 0 \mod 4$,}\\
& \text{or $4n - 5 \le i \le 8n - 13$ and $i \equiv -1 \mod 4$;}\\
\ZM/2 & \text{if $2n - 4 < i < 4n - 6$ and $i \equiv 2 \mod 4$;}\\
& \text{or $4n - 6 < i < 6n - 8$ and $i \equiv 2 \mod 4$;}\\
(\ZM/2)^2 & \text{if $i = 4n - 6$ and $n$ is even;}\\
\ZM/4 & \text{if $i = 4n - 6$ and $n$ is odd;}\\
0 & \text{otherwise.}
\end{cases}
\oplus
\ZM\ 
\begin{array}{ll}
\text{if $i = 2n - 4$}\\
\text{or $i =  6n - 7$.}
\end{array}
\]

\subsection{Type $E_6$}

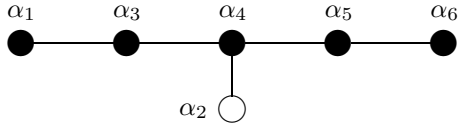
\begin{center}
\begin{picture}(200,45)
\put( 40, 30){\circle*{10}}
\put( 45, 30){\line(1,0){30}}
\put( 80, 30){\circle*{10}}
\put( 85, 30){\line(1,0){30}}
\put(120, 30){\circle*{10}}
\put(125, 30){\line(1,0){30}}
\put(160, 30){\circle*{10}}
\put(165, 30){\line(1,0){30}}
\put(200, 30){\circle*{10}}
\put(120,  5){\circle{10}}
\put(120, 25){\line(0,-1){15}}
\put( 35, 40){$\a_1$}
\put( 75, 40){$\a_3$}
\put(115, 40){$\a_4$}
\put(155, 40){$\a_5$}
\put(195, 40){$\a_6$}
\put(100,  3){$\a_2$}
\end{picture}
\end{center}

We have $h = h^\vee = 12$, and $d = 22$. The Cartan matrix has cokernel
isomorphic to $\ZM/3$.

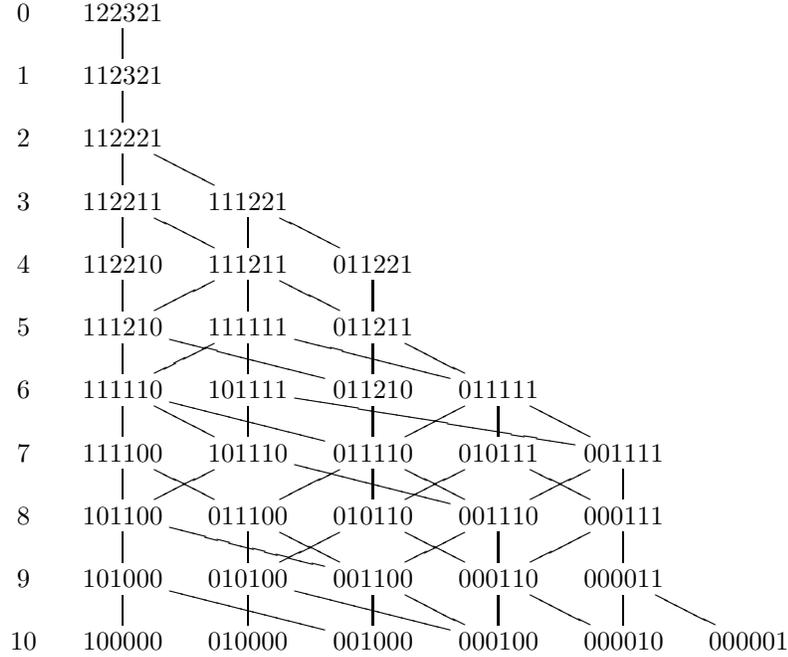
\[
\xymatrix @=.4cm{
0 &
122321\ar@{-}[d]\\
1 &
112321\ar@{-}[d]\\
2 &
112221\ar@{-}[d]\ar@{-}[dr]\\
3 &
112211\ar@{-}[d]\ar@{-}[dr]&
111221\ar@{-}[d]\ar@{-}[dr]\\
4 &
112210\ar@{-}[d]&
111211\ar@{-}[dl]\ar@{-}[d]\ar@{-}[dr]&
011221\ar@{-}[d]\\
5 &
111210\ar@{-}[d]\ar@{-}[drr]&
111111\ar@{-}[dl]\ar@{-}[d]\ar@{-}[drr]&
011211\ar@{-}[d]\ar@{-}[dr]\\
6 &
111110\ar@{-}[d]\ar@{-}[dr]\ar@{-}[drr]&
101111\ar@{-}[d]\ar@{-}[drrr]&
011210\ar@{-}[d]&
011111\ar@{-}[dl]\ar@{-}[d]\ar@{-}[dr]\\
7 &
111100\ar@{-}[d]\ar@{-}[dr]&
101110\ar@{-}[dl]\ar@{-}[drr]&
011110\ar@{-}[dl]\ar@{-}[d]\ar@{-}[dr]&
010111\ar@{-}[dl]\ar@{-}[dr]&
001111\ar@{-}[dl]\ar@{-}[d]\\
8 &
101100\ar@{-}[d]\ar@{-}[drr]&
011100\ar@{-}[d]\ar@{-}[dr]&
010110\ar@{-}[dl]\ar@{-}[dr]&
001110\ar@{-}[dl]\ar@{-}[d]&
000111\ar@{-}[dl]\ar@{-}[d]\\
9 &
101000\ar@{-}[d]\ar@{-}[drr]&
010100\ar@{-}[d]\ar@{-}[drr]&
001100\ar@{-}[d]\ar@{-}[dr]&
000110\ar@{-}[d]\ar@{-}[dr]&
000011\ar@{-}[d]\ar@{-}[dr]\\
10&
100000&
010000&
001000&
000100&
000010&
000001
}
\]

\[
\DC_1 = 
\DC_2 = \begin{pmatrix}  1 \end{pmatrix}\quad
\DC_3 = \begin{pmatrix}  1 \\
     1 \end{pmatrix}\quad
\DC_4 = \begin{pmatrix}  1 & 0 \\
     1 & 1 \\
     0 & 1 \end{pmatrix}\quad
\DC_5 = \begin{pmatrix}  1 & 1 & 0 \\
     0 & 1 & 0 \\
     0 & 1 & 1 \end{pmatrix}\] \[
\DC_6 = \begin{pmatrix}  1 & 1 & 0 \\
     0 & 1 & 0 \\
     1 & 0 & 1 \\
     0 & 1 & 1 \end{pmatrix}\quad
\DC_7 = \begin{pmatrix}  1 & 0 & 0 & 0 \\
     1 & 1 & 0 & 0 \\
     1 & 0 & 1 & 1 \\
     0 & 0 & 0 & 1 \\
     0 & 1 & 0 & 1 \end{pmatrix}\quad
\DC_8 = \begin{pmatrix}  1 & 1 & 0 & 0 & 0 \\
     1 & 0 & 1 & 0 & 0 \\
     0 & 0 & 1 & 1 & 0 \\
     0 & 1 & 1 & 0 & 1 \\
     0 & 0 & 0 & 1 & 1 \end{pmatrix}\] \[
\DC_9 = \begin{pmatrix}  1 & 0 & 0 & 0 & 0 \\
     0 & 1 & 1 & 0 & 0 \\
     1 & 1 & 0 & 1 & 0 \\
     0 & 0 & 1 & 1 & 1 \\
     0 & 0 & 0 & 0 & 1 \end{pmatrix}\quad
\DC_{10} = \begin{pmatrix}  1 & 0 & 0 & 0 & 0 \\
     0 & 1 & 0 & 0 & 0 \\
     1 & 0 & 1 & 0 & 0 \\
     0 & 1 & 1 & 1 & 0 \\
     0 & 0 & 0 & 1 & 1 \\
     0 & 0 & 0 & 0 & 1 \end{pmatrix}
\]

\[
H^i(\OC_\mini, \ZM) \simeq
\begin{cases}
\ZM & \text{for $i = 0,\ 6,\ 8,\ 12,\ 14,\ 20,\ 23,\ 29,\ 31,\ 35,\ 37,\ 43$}\\
\ZM/3 & \text{for $i = 16,\ 22,\ 28$}\\
\ZM/2 & \text{for $i = 18,\ 26$}\\
0 & \text{otherwise}
\end{cases}
\]

\subsection{Type $E_7$}

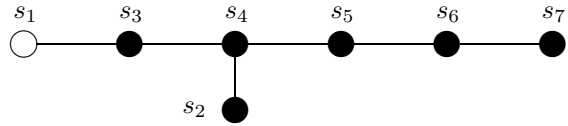
\begin{center}
\begin{picture}(240,45)
\put( 40, 30){\circle{10}}
\put( 45, 30){\line(1,0){30}}
\put( 80, 30){\circle*{10}}
\put( 85, 30){\line(1,0){30}}
\put(120, 30){\circle*{10}}
\put(125, 30){\line(1,0){30}}
\put(160, 30){\circle*{10}}
\put(165, 30){\line(1,0){30}}
\put(200, 30){\circle*{10}}
\put(205, 30){\line(1,0){30}}
\put(240, 30){\circle*{10}}
\put(120, 5){\circle*{10}}
\put(120, 25){\line(0,-1){15}}
\put( 36, 40){$s_1$}
\put( 76, 40){$s_3$}
\put(116, 40){$s_4$}
\put(156, 40){$s_5$}
\put(196, 40){$s_6$}
\put(236, 40){$s_7$}
\put(100, 4){$s_2$}
\end{picture}
\end{center}

We have $h = h^\vee = 18$, and $d = 34$. The Cartan matrix has cokernel
isomorphic to $\ZM/2$.

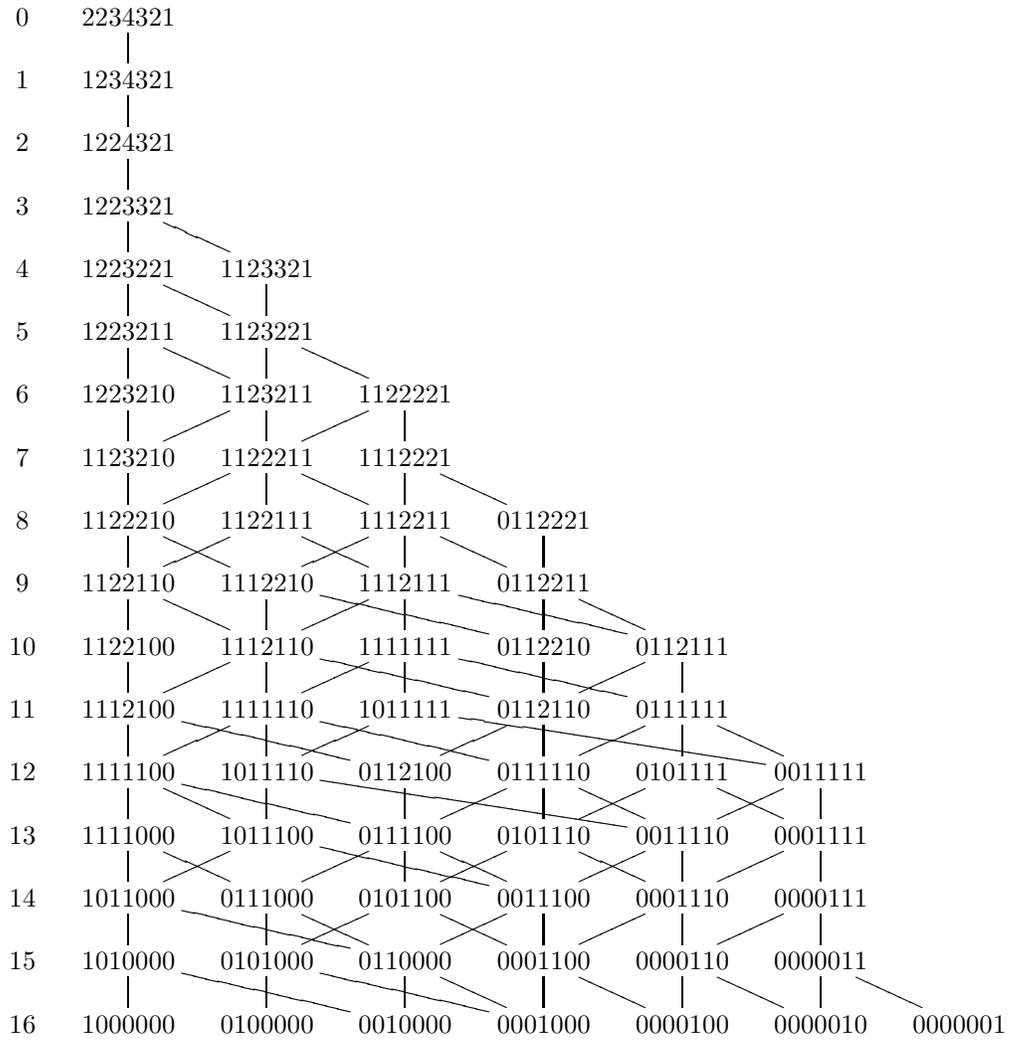
\[
\xymatrix @=.4cm{
0 &
2234321\ar@{-}[d]\\
1 &
1234321\ar@{-}[d]\\
2 &
1224321\ar@{-}[d]\\
3 &
1223321\ar@{-}[d]\ar@{-}[dr]\\
4 &
1223221\ar@{-}[d]\ar@{-}[dr]&
1123321\ar@{-}[d]\\
5 &
1223211\ar@{-}[d]\ar@{-}[dr]&
1123221\ar@{-}[d]\ar@{-}[dr]\\
6 &
1223210\ar@{-}[d]&
1123211\ar@{-}[dl]\ar@{-}[d]&
1122221\ar@{-}[dl]\ar@{-}[d]\\
7 &
1123210\ar@{-}[d]&
1122211\ar@{-}[dl]\ar@{-}[d]\ar@{-}[dr]&
1112221\ar@{-}[d]\ar@{-}[dr]\\
8 &
1122210\ar@{-}[d]\ar@{-}[dr]&
1122111\ar@{-}[dl]\ar@{-}[dr]&
1112211\ar@{-}[dl]\ar@{-}[d]\ar@{-}[dr]&
0112221\ar@{-}[d]\\
9 &
1122110\ar@{-}[d]\ar@{-}[dr]&
1112210\ar@{-}[d]\ar@{-}[drr]&
1112111\ar@{-}[dl]\ar@{-}[d]\ar@{-}[drr]&
0112211\ar@{-}[d]\ar@{-}[dr]\\
10&
1122100\ar@{-}[d]&
1112110\ar@{-}[dl]\ar@{-}[d]\ar@{-}[drr]&
1111111\ar@{-}[dl]\ar@{-}[d]\ar@{-}[drr]&
0112210\ar@{-}[d]&
0112111\ar@{-}[dl]\ar@{-}[d]\\
11&
1112100\ar@{-}[d]\ar@{-}[drr]&
1111110\ar@{-}[dl]\ar@{-}[d]\ar@{-}[drr]&
1011111\ar@{-}[dl]\ar@{-}[drrr]&
0112110\ar@{-}[dl]\ar@{-}[d]&
0111111\ar@{-}[dl]\ar@{-}[d]\ar@{-}[dr]\\
12&
1111100\ar@{-}[d]\ar@{-}[dr]\ar@{-}[drr]&
1011110\ar@{-}[d]\ar@{-}[drrr]&
0112100\ar@{-}[d]&
0111110\ar@{-}[dl]\ar@{-}[d]\ar@{-}[dr]&
0101111\ar@{-}[dl]\ar@{-}[dr]&
0011111\ar@{-}[dl]\ar@{-}[d]\\
13&
1111000\ar@{-}[d]\ar@{-}[dr]&
1011100\ar@{-}[dl]\ar@{-}[drr]&
0111100\ar@{-}[dl]\ar@{-}[d]\ar@{-}[dr]&
0101110\ar@{-}[dl]\ar@{-}[dr]&
0011110\ar@{-}[dl]\ar@{-}[d]&
0001111\ar@{-}[dl]\ar@{-}[d]\\
14&
1011000\ar@{-}[d]\ar@{-}[drr]&
0111000\ar@{-}[d]\ar@{-}[dr]&
0101100\ar@{-}[dl]\ar@{-}[dr]&
0011100\ar@{-}[dl]\ar@{-}[d]&
0001110\ar@{-}[dl]\ar@{-}[d]&
0000111\ar@{-}[dl]\ar@{-}[d]\\
15&
1010000\ar@{-}[d]\ar@{-}[drr]&
0101000\ar@{-}[d]\ar@{-}[drr]&
0110000\ar@{-}[d]\ar@{-}[dr]&
0001100\ar@{-}[d]\ar@{-}[dr]&
0000110\ar@{-}[d]\ar@{-}[dr]&
0000011\ar@{-}[d]\ar@{-}[dr]\\
16&
1000000&
0100000&
0010000&
0001000&
0000100&
0000010&
0000001
}
\]

\newpage

\[ \DC_1 = \DC_2 = \DC_3 = \begin{pmatrix}  1 \end{pmatrix}\quad
\DC_4 = \begin{pmatrix}  1 \\
     1 \end{pmatrix}\quad
\DC_5 = \begin{pmatrix}  1 & 0 \\
     1 & 1 \end{pmatrix}\quad
\DC_6 = \begin{pmatrix}  1 & 0 \\
     1 & 1 \\
     0 & 1 \end{pmatrix}\] \[
\DC_7 = \begin{pmatrix}  1 & 1 & 0 \\
     0 & 1 & 1 \\
     0 & 0 & 1 \end{pmatrix}\quad
\DC_8 = \begin{pmatrix}  1 & 1 & 0 \\
     0 & 1 & 0 \\
     0 & 1 & 1 \\
     0 & 0 & 1 \end{pmatrix}\quad
\DC_9 = \begin{pmatrix}  1 & 1 & 0 & 0 \\
     1 & 0 & 1 & 0 \\
     0 & 1 & 1 & 0 \\
     0 & 0 & 1 & 1 \end{pmatrix}\quad
\DC_{10} = \begin{pmatrix}  1 & 0 & 0 & 0 \\
     1 & 1 & 1 & 0 \\
     0 & 0 & 1 & 0 \\
     0 & 1 & 0 & 1 \\
     0 & 0 & 1 & 1 \end{pmatrix}\] \[
\DC_{11} = \begin{pmatrix}  1 & 1 & 0 & 0 & 0 \\
     0 & 1 & 1 & 0 & 0 \\
     0 & 0 & 1 & 0 & 0 \\
     0 & 1 & 0 & 1 & 1 \\
     0 & 0 & 1 & 0 & 1 \end{pmatrix}\quad
\DC_{12} = \begin{pmatrix}  1 & 1 & 0 & 0 & 0 \\
     0 & 1 & 1 & 0 & 0 \\
     1 & 0 & 0 & 1 & 0 \\
     0 & 1 & 0 & 1 & 1 \\
     0 & 0 & 0 & 0 & 1 \\
     0 & 0 & 1 & 0 & 1 \end{pmatrix}\] \[
\DC_{13} = \begin{pmatrix}  1 & 0 & 0 & 0 & 0 & 0 \\
     1 & 1 & 0 & 0 & 0 & 0 \\
     1 & 0 & 1 & 1 & 0 & 0 \\
     0 & 0 & 0 & 1 & 1 & 0 \\
     0 & 1 & 0 & 1 & 0 & 1 \\
     0 & 0 & 0 & 0 & 1 & 1 \end{pmatrix}\quad
\DC_{14} = \begin{pmatrix}  1 & 1 & 0 & 0 & 0 & 0 \\
     1 & 0 & 1 & 0 & 0 & 0 \\
     0 & 0 & 1 & 1 & 0 & 0 \\
     0 & 1 & 1 & 0 & 1 & 0 \\
     0 & 0 & 0 & 1 & 1 & 1 \\
     0 & 0 & 0 & 0 & 0 & 1 \end{pmatrix}\] \[
\DC_{15} = \begin{pmatrix}  1 & 0 & 0 & 0 & 0 & 0 \\
     0 & 1 & 1 & 0 & 0 & 0 \\
     1 & 1 & 0 & 1 & 0 & 0 \\
     0 & 0 & 1 & 1 & 1 & 0 \\
     0 & 0 & 0 & 0 & 1 & 1 \\
     0 & 0 & 0 & 0 & 0 & 1 \end{pmatrix}\quad
\DC_{16} = \begin{pmatrix}  1 & 0 & 0 & 0 & 0 & 0 \\
     0 & 1 & 0 & 0 & 0 & 0 \\
     1 & 0 & 1 & 0 & 0 & 0 \\
     0 & 1 & 1 & 1 & 0 & 0 \\
     0 & 0 & 0 & 1 & 1 & 0 \\
     0 & 0 & 0 & 0 & 1 & 1 \\
     0 & 0 & 0 & 0 & 0 & 1 \end{pmatrix}
\]

\[
H^i(\OC_\mini, \ZM) \simeq
\begin{cases}
\ZM & \text{for } i = 0,\ 8,\ 12,\ 16,\ 20,\ 24,\ 32,\\ 
    & \qquad \quad 35,\ 43,\ 47,\ 51,\ 55,\ 59,\ 67\\
\ZM/2 & \text{for } i = 18,\ 26,\ 30,\ 34,\ 38,\ 42,\ 50\\
\ZM/3 & \text{for } i = 28,\ 40\\
0 & \text{otherwise}
\end{cases}
\]

\subsection{Type $E_8$}

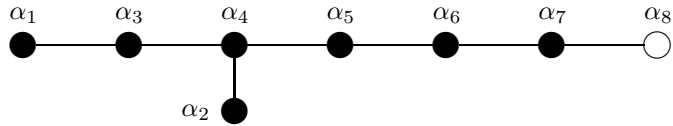
\begin{center}
\begin{picture}(280,45)
\put( 40, 30){\circle*{10}}
\put( 45, 30){\line(1,0){30}}
\put( 80, 30){\circle*{10}}
\put( 85, 30){\line(1,0){30}}
\put(120, 30){\circle*{10}}
\put(125, 30){\line(1,0){30}}
\put(160, 30){\circle*{10}}
\put(165, 30){\line(1,0){30}}
\put(200, 30){\circle*{10}}
\put(205, 30){\line(1,0){30}}
\put(240, 30){\circle*{10}}
\put(245, 30){\line(1,0){30}}
\put(280, 30){\circle{10}}
\put(120,  5){\circle*{10}}
\put(120, 25){\line(0,-1){15}}
\put( 35, 40){$\a_1$}
\put( 75, 40){$\a_3$}
\put(115, 40){$\a_4$}
\put(155, 40){$\a_5$}
\put(195, 40){$\a_6$}
\put(235, 40){$\a_7$}
\put(275, 40){$\a_8$}
\put(100,  3){$\a_2$}
\end{picture}
\end{center}

We have $h = h^\vee = 30$, and $d = 58$. The Cartan matrix is an
isomorphism.

\newpage

{\scriptsize
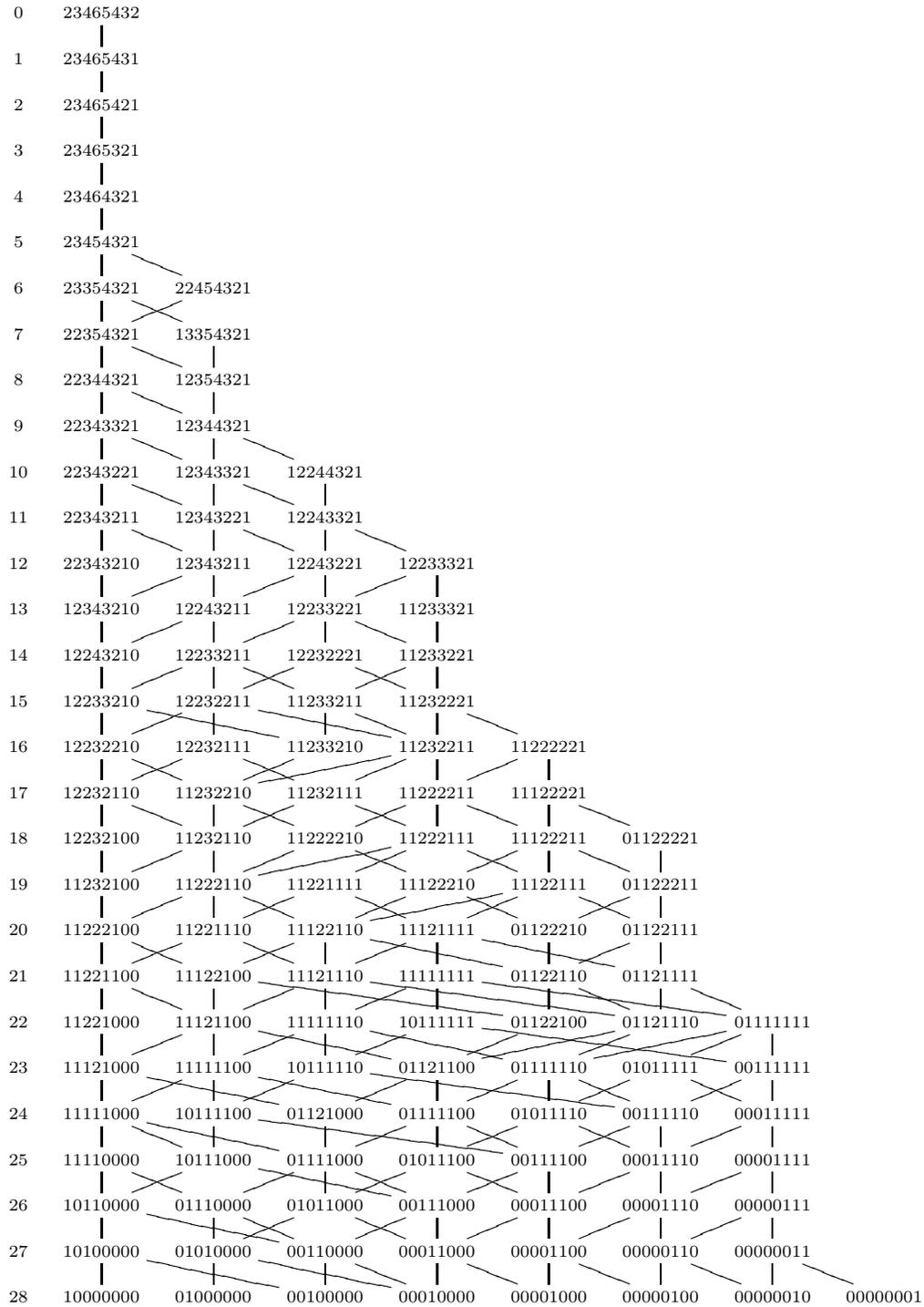
\[
\xymatrix @=.3cm{
0 &
23465432\ar@{-}[d]\\
1 &
23465431\ar@{-}[d]\\
2 &
23465421\ar@{-}[d]\\
3 &
23465321\ar@{-}[d]\\
4 &
23464321\ar@{-}[d]\\
5 &
23454321\ar@{-}[d]\ar@{-}[dr]\\
6 &
23354321\ar@{-}[d]\ar@{-}[dr]&
22454321\ar@{-}[dl]\\
7 &
22354321\ar@{-}[d]\ar@{-}[dr]&
13354321\ar@{-}[d]\\
8 &
22344321\ar@{-}[d]\ar@{-}[dr]&
12354321\ar@{-}[d]\\
9 &
22343321\ar@{-}[d]\ar@{-}[dr]&
12344321\ar@{-}[d]\ar@{-}[dr]\\
10&
22343221\ar@{-}[d]\ar@{-}[dr]&
12343321\ar@{-}[d]\ar@{-}[dr]&
12244321\ar@{-}[d]\\
11&
22343211\ar@{-}[d]\ar@{-}[dr]&
12343221\ar@{-}[d]\ar@{-}[dr]&
12243321\ar@{-}[d]\ar@{-}[dr]\\
12&
22343210\ar@{-}[d]&
12343211\ar@{-}[dl]\ar@{-}[d]&
12243221\ar@{-}[dl]\ar@{-}[d]&
12233321\ar@{-}[dl]\ar@{-}[d]\\
13&
12343210\ar@{-}[d]&
12243211\ar@{-}[dl]\ar@{-}[d]&
12233221\ar@{-}[dl]\ar@{-}[d]\ar@{-}[dr]&
11233321\ar@{-}[d]\\
14&
12243210\ar@{-}[d]&
12233211\ar@{-}[dl]\ar@{-}[d]\ar@{-}[dr]&
12232221\ar@{-}[dl]\ar@{-}[dr]&
11233221\ar@{-}[dl]\ar@{-}[d]\\
15&
12233210\ar@{-}[d]\ar@{-}[drr]&
12232211\ar@{-}[dl]\ar@{-}[d]\ar@{-}[drr]&
11233211\ar@{-}[d]\ar@{-}[dr]&
11232221\ar@{-}[d]\ar@{-}[dr]\\
16&
12232210\ar@{-}[d]\ar@{-}[dr]&
12232111\ar@{-}[dl]\ar@{-}[dr]&
11233210\ar@{-}[dl]&
11232211\ar@{-}[dll]\ar@{-}[dl]\ar@{-}[d]&
11222221\ar@{-}[dl]\ar@{-}[d]\\
17&
12232110\ar@{-}[d]\ar@{-}[dr]&
11232210\ar@{-}[d]\ar@{-}[dr]&
11232111\ar@{-}[dl]\ar@{-}[dr]&
11222211\ar@{-}[dl]\ar@{-}[d]\ar@{-}[dr]&
11122221\ar@{-}[d]\ar@{-}[dr]\\
18&
12232100\ar@{-}[d]&
11232110\ar@{-}[dl]\ar@{-}[d]&
11222210\ar@{-}[dl]\ar@{-}[dr]&
11222111\ar@{-}[dll]\ar@{-}[dl]\ar@{-}[dr]&
11122211\ar@{-}[dl]\ar@{-}[d]\ar@{-}[dr]&
01122221\ar@{-}[d]\\
19&
11232100\ar@{-}[d]&
11222110\ar@{-}[dl]\ar@{-}[d]\ar@{-}[dr]&
11221111\ar@{-}[dl]\ar@{-}[dr]&
11122210\ar@{-}[dl]\ar@{-}[dr]&
11122111\ar@{-}[dll]\ar@{-}[dl]\ar@{-}[dr]&
01122211\ar@{-}[dl]\ar@{-}[d]\\
20&
11222100\ar@{-}[d]\ar@{-}[dr]&
11221110\ar@{-}[dl]\ar@{-}[dr]&
11122110\ar@{-}[dl]\ar@{-}[d]\ar@{-}[drr]&
11121111\ar@{-}[dl]\ar@{-}[d]\ar@{-}[drr]&
01122210\ar@{-}[d]&
01122111\ar@{-}[dl]\ar@{-}[d]\\
21&
11221100\ar@{-}[d]\ar@{-}[dr]&
11122100\ar@{-}[d]\ar@{-}[drrr]&
11121110\ar@{-}[dl]\ar@{-}[d]\ar@{-}[drrr]&
11111111\ar@{-}[dl]\ar@{-}[d]\ar@{-}[drrr]&
01122110\ar@{-}[d]\ar@{-}[dr]&
01121111\ar@{-}[d]\ar@{-}[dr]\\
22&
11221000\ar@{-}[d]&
11121100\ar@{-}[dl]\ar@{-}[d]\ar@{-}[drr]&
11111110\ar@{-}[dl]\ar@{-}[d]\ar@{-}[drr]&
10111111\ar@{-}[dl]\ar@{-}[drrr]&
01122100\ar@{-}[dl]&
01121110\ar@{-}[dll]\ar@{-}[dl]&
01111111\ar@{-}[dll]\ar@{-}[dl]\ar@{-}[d]\\
23&
11121000\ar@{-}[d]\ar@{-}[drr]&
11111100\ar@{-}[dl]\ar@{-}[d]\ar@{-}[drr]&
10111110\ar@{-}[dl]\ar@{-}[drrr]&
01121100\ar@{-}[dl]\ar@{-}[d]&
01111110\ar@{-}[dl]\ar@{-}[d]\ar@{-}[dr]&
01011111\ar@{-}[dl]\ar@{-}[dr]&
00111111\ar@{-}[dl]\ar@{-}[d]\\
24&
11111000\ar@{-}[d]\ar@{-}[dr]\ar@{-}[drr]&
10111100\ar@{-}[d]\ar@{-}[drrr]&
01121000\ar@{-}[d]&
01111100\ar@{-}[dl]\ar@{-}[d]\ar@{-}[dr]&
01011110\ar@{-}[dl]\ar@{-}[dr]&
00111110\ar@{-}[dl]\ar@{-}[d]&
00011111\ar@{-}[dl]\ar@{-}[d]\\
25&
11110000\ar@{-}[d]\ar@{-}[dr]&
10111000\ar@{-}[dl]\ar@{-}[drr]&
01111000\ar@{-}[dl]\ar@{-}[d]\ar@{-}[dr]&
01011100\ar@{-}[dl]\ar@{-}[dr]&
00111100\ar@{-}[dl]\ar@{-}[d]&
00011110\ar@{-}[dl]\ar@{-}[d]&
00001111\ar@{-}[dl]\ar@{-}[d]\\
26&
10110000\ar@{-}[d]\ar@{-}[drr]&
01110000\ar@{-}[d]\ar@{-}[dr]&
01011000\ar@{-}[dl]\ar@{-}[dr]&
00111000\ar@{-}[dl]\ar@{-}[d]&
00011100\ar@{-}[dl]\ar@{-}[d]&
00001110\ar@{-}[dl]\ar@{-}[d]&
00000111\ar@{-}[dl]\ar@{-}[d]\\
27&
10100000\ar@{-}[d]\ar@{-}[drr]&
01010000\ar@{-}[d]\ar@{-}[drr]&
00110000\ar@{-}[d]\ar@{-}[dr]&
00011000\ar@{-}[d]\ar@{-}[dr]&
00001100\ar@{-}[d]\ar@{-}[dr]&
00000110\ar@{-}[d]\ar@{-}[dr]&
00000011\ar@{-}[d]\ar@{-}[dr]\\
28&
10000000&
01000000&
00100000&
00010000&
00001000&
00000100&
00000010&
00000001
}
\]
}

\newpage

{\scriptsize
\[
\DC_1 = 
\DC_2 = 
\DC_3 = 
\DC_4 = 
\DC_5 = \begin{pmatrix}  1 \end{pmatrix}\quad
\DC_6 = \begin{pmatrix}  1 \\
     1 \end{pmatrix}\quad
\DC_7 = \begin{pmatrix}  1 & 1 \\
     1 & 0 \end{pmatrix}\quad
\DC_8 = 
\DC_9 = \begin{pmatrix}  1 & 0 \\
     1 & 1 \end{pmatrix}\] \[
\DC_{10} = \begin{pmatrix}  1 & 0 \\
     1 & 1 \\
     0 & 1 \end{pmatrix}\quad
\DC_{11} = \begin{pmatrix}  1 & 0 & 0 \\
     1 & 1 & 0 \\
     0 & 1 & 1 \end{pmatrix}\quad
\DC_{12} = \begin{pmatrix}  1 & 0 & 0 \\
     1 & 1 & 0 \\
     0 & 1 & 1 \\
     0 & 0 & 1 \end{pmatrix}\quad
\DC_{13} = \begin{pmatrix}  1 & 1 & 0 & 0 \\
     0 & 1 & 1 & 0 \\
     0 & 0 & 1 & 1 \\
     0 & 0 & 0 & 1 \end{pmatrix}\] \[
\DC_{14} = \begin{pmatrix}  1 & 1 & 0 & 0 \\
     0 & 1 & 1 & 0 \\
     0 & 0 & 1 & 0 \\
     0 & 0 & 1 & 1 \end{pmatrix}\quad
\DC_{15} = \begin{pmatrix}  1 & 1 & 0 & 0 \\
     0 & 1 & 1 & 0 \\
     0 & 1 & 0 & 1 \\
     0 & 0 & 1 & 1 \end{pmatrix}\quad
\DC_{16} = \begin{pmatrix}  1 & 1 & 0 & 0 \\
     0 & 1 & 0 & 0 \\
     1 & 0 & 1 & 0 \\
     0 & 1 & 1 & 1 \\
     0 & 0 & 0 & 1 \end{pmatrix}\] \[
\DC_{17} = \begin{pmatrix}  1 & 1 & 0 & 0 & 0 \\
     1 & 0 & 1 & 1 & 0 \\
     0 & 1 & 0 & 1 & 0 \\
     0 & 0 & 0 & 1 & 1 \\
     0 & 0 & 0 & 0 & 1 \end{pmatrix}\quad
\DC_{18} = \begin{pmatrix}  1 & 0 & 0 & 0 & 0 \\
     1 & 1 & 1 & 0 & 0 \\
     0 & 1 & 0 & 1 & 0 \\
     0 & 0 & 1 & 1 & 0 \\
     0 & 0 & 0 & 1 & 1 \\
     0 & 0 & 0 & 0 & 1 \end{pmatrix}\quad
\DC_{19} = \begin{pmatrix}  1 & 1 & 0 & 0 & 0 & 0 \\
     0 & 1 & 1 & 1 & 0 & 0 \\
     0 & 0 & 0 & 1 & 0 & 0 \\
     0 & 0 & 1 & 0 & 1 & 0 \\
     0 & 0 & 0 & 1 & 1 & 0 \\
     0 & 0 & 0 & 0 & 1 & 1 \end{pmatrix}\] \[
\DC_{20} = \begin{pmatrix}  1 & 1 & 0 & 0 & 0 & 0 \\
     0 & 1 & 1 & 0 & 0 & 0 \\
     0 & 1 & 0 & 1 & 1 & 0 \\
     0 & 0 & 1 & 0 & 1 & 0 \\
     0 & 0 & 0 & 1 & 0 & 1 \\
     0 & 0 & 0 & 0 & 1 & 1 \end{pmatrix}\quad
\DC_{21} = \begin{pmatrix}  1 & 1 & 0 & 0 & 0 & 0 \\
     1 & 0 & 1 & 0 & 0 & 0 \\
     0 & 1 & 1 & 1 & 0 & 0 \\
     0 & 0 & 0 & 1 & 0 & 0 \\
     0 & 0 & 1 & 0 & 1 & 1 \\
     0 & 0 & 0 & 1 & 0 & 1 \end{pmatrix}\quad
\DC_{22} = \begin{pmatrix}  1 & 0 & 0 & 0 & 0 & 0 \\
     1 & 1 & 1 & 0 & 0 & 0 \\
     0 & 0 & 1 & 1 & 0 & 0 \\
     0 & 0 & 0 & 1 & 0 & 0 \\
     0 & 1 & 0 & 0 & 1 & 0 \\
     0 & 0 & 1 & 0 & 1 & 1 \\
     0 & 0 & 0 & 1 & 0 & 1 \end{pmatrix}\] \[
\DC_{23} = \begin{pmatrix}  1 & 1 & 0 & 0 & 0 & 0 & 0 \\
     0 & 1 & 1 & 0 & 0 & 0 & 0 \\
     0 & 0 & 1 & 1 & 0 & 0 & 0 \\
     0 & 1 & 0 & 0 & 1 & 1 & 0 \\
     0 & 0 & 1 & 0 & 0 & 1 & 1 \\
     0 & 0 & 0 & 0 & 0 & 0 & 1 \\
     0 & 0 & 0 & 1 & 0 & 0 & 1 \end{pmatrix}\quad
\DC_{24} = \begin{pmatrix}  1 & 1 & 0 & 0 & 0 & 0 & 0 \\
     0 & 1 & 1 & 0 & 0 & 0 & 0 \\
     1 & 0 & 0 & 1 & 0 & 0 & 0 \\
     0 & 1 & 0 & 1 & 1 & 0 & 0 \\
     0 & 0 & 0 & 0 & 1 & 1 & 0 \\
     0 & 0 & 1 & 0 & 1 & 0 & 1 \\
     0 & 0 & 0 & 0 & 0 & 1 & 1 \end{pmatrix}\] \[
\DC_{25} = \begin{pmatrix}  1 & 0 & 0 & 0 & 0 & 0 & 0 \\
     1 & 1 & 0 & 0 & 0 & 0 & 0 \\
     1 & 0 & 1 & 1 & 0 & 0 & 0 \\
     0 & 0 & 0 & 1 & 1 & 0 & 0 \\
     0 & 1 & 0 & 1 & 0 & 1 & 0 \\
     0 & 0 & 0 & 0 & 1 & 1 & 1 \\
     0 & 0 & 0 & 0 & 0 & 0 & 1 \end{pmatrix}\quad
\DC_{26} = \begin{pmatrix}  1 & 1 & 0 & 0 & 0 & 0 & 0 \\
     1 & 0 & 1 & 0 & 0 & 0 & 0 \\
     0 & 0 & 1 & 1 & 0 & 0 & 0 \\
     0 & 1 & 1 & 0 & 1 & 0 & 0 \\
     0 & 0 & 0 & 1 & 1 & 1 & 0 \\
     0 & 0 & 0 & 0 & 0 & 1 & 1 \\
     0 & 0 & 0 & 0 & 0 & 0 & 1 \end{pmatrix}\] \[
\DC_{27} = \begin{pmatrix}  1 & 0 & 0 & 0 & 0 & 0 & 0 \\
     0 & 1 & 1 & 0 & 0 & 0 & 0 \\
     1 & 1 & 0 & 1 & 0 & 0 & 0 \\
     0 & 0 & 1 & 1 & 1 & 0 & 0 \\
     0 & 0 & 0 & 0 & 1 & 1 & 0 \\
     0 & 0 & 0 & 0 & 0 & 1 & 1 \\
     0 & 0 & 0 & 0 & 0 & 0 & 1 \end{pmatrix}\quad
\DC_{28} = \begin{pmatrix}  1 & 0 & 0 & 0 & 0 & 0 & 0 \\
     0 & 1 & 0 & 0 & 0 & 0 & 0 \\
     1 & 0 & 1 & 0 & 0 & 0 & 0 \\
     0 & 1 & 1 & 1 & 0 & 0 & 0 \\
     0 & 0 & 0 & 1 & 1 & 0 & 0 \\
     0 & 0 & 0 & 0 & 1 & 1 & 0 \\
     0 & 0 & 0 & 0 & 0 & 1 & 1 \\
     0 & 0 & 0 & 0 & 0 & 0 & 1 \end{pmatrix}\quad
\]
}

\[
H^i(\OC_\mini, \ZM) \simeq
\begin{cases}
\ZM & \text{for } i = 0,\ 12,\ 20,\ 24,\ 32,\ 36,\ 44,\ 56,\\ 
    & \qquad \quad 59,\ 71,\ 79,\ 83,\ 91,\ 95,\ 103, 115\\
\ZM/2 & \text{for } i = 30,\ 42,\ 50,\ 54,\ 62,\ 66,\ 74,\ 86\\
\ZM/3 & \text{for } i = 40,\ 52,\ 64,\ 76\\
\ZM/5 & \text{for } i = 48,\ 68\\
0 & \text{otherwise}
\end{cases}
\]

\newpage

\subsection{Type $F_4$}

\begin{center}
\begin{picture}(200,30)(-10,0)
\put( 40, 10){\circle{10}}
\put( 45, 10){\line(1,0){30}}
\put( 80, 10){\circle*{10}}
\put( 84, 11.5){\line(1,0){32}}
\put( 84,  8.5){\line(1,0){32}}
\put( 95,  5.5){\LARGE{$>$}}
\put(120, 10){\circle*{10}}
\put(125, 10){\line(1,0){30}}
\put(160, 10){\circle*{10}}
\put( 35, 20){$\a_1$}
\put( 75, 20){$\a_2$}
\put(115, 20){$\a_3$}
\put(155, 20){$\a_4$}
\end{picture}
\end{center}

We have $h = 12$, $h^\vee = 9$ and $d = 16$.

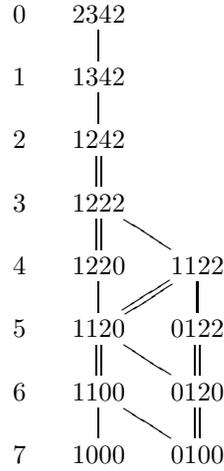
\[
\xymatrix @=.4cm{
0&
2342\ar@{-}[d]\\
1&
1342\ar@{-}[d]\\
2&
1242\ar@{=}[d]\\
3&
1222\ar@{=}[d]\ar@{-}[dr]\\
4&
1220\ar@{-}[d]&
1122\ar@{=}[dl]\ar@{-}[d]\\
5&
1120\ar@{=}[d]\ar@{-}[dr]&
0122\ar@{=}[d]\\
6&
1100\ar@{-}[d]\ar@{-}[dr]&
0120\ar@{=}[d]\\
7&
1000&
0100
}
\]

We have
\[
\DC_1 = \DC_2 = (1)\quad \DC_3 = (2)
\quad \DC_4 = \begin{pmatrix}
2\\
1
\end{pmatrix}
\quad \DC_5 = \begin{pmatrix}
1&2\\
0&1
\end{pmatrix}
\quad \DC_6 = \begin{pmatrix}
2&0\\
1&2
\end{pmatrix}
\quad \DC_7 = \begin{pmatrix}
1&0\\
1&2
\end{pmatrix}
\]

The type of $\Phi'$ is $A_2$, so we have
\[
\DC_8 =
\begin{pmatrix}
2&1\\
1&2
\end{pmatrix}
\]

The matrices of the last differentials are transposed to the first ones.

\[
H^i(\OC_\mini, \ZM) \simeq
\begin{cases}
\ZM & \text{for } i = 0,\ 8,\ 23,\ 31\\
\ZM/2 & \text{for } i = 6,\ 14,\ 18,\ 26\\
\ZM/4 & \text{for } i = 12,\ 20\\
\ZM/3 & \text{for } i = 16\\
0 & \text{otherwise}
\end{cases}
\]

\subsection{Type $G_2$}

\begin{center}
\begin{picture}(100,30)
\put( 40, 10){\circle{10}}
\put(44.5,12){\line(1,0){32}}
\put( 45, 10){\line(1,0){30}}
\put(44.5, 8){\line(1,0){32}}
\put( 55,5.5){\LARGE{$>$}}
\put( 80, 10){\circle*{10}}
\put( 35, 20){$\a_1$}
\put( 75, 20){$\a_2$}
\end{picture}
\end{center}

We have $h = 6$, $h^\vee = 4$, and $d = 6$.
The root system $\Phi'$ is of type $A_1$. Its Cartan matrix
has cokernel $\ZM/2$.

$$
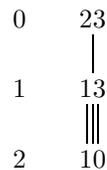
\xymatrix @=.5cm{
0&
23\ar@{-}[d]\\
1&
13\ar@3{-}[d]\\
2&
10
}
$$

We have
\[
\DC_1 = (1) \quad 
\DC_2 = (3) \quad 
\DC_3 = (2) \quad 
\DC_4 = (3) \quad 
\DC_5 = (1)
\]

\[
H^i(\OC_\mini, \ZM) \simeq
\begin{cases}
\ZM & \text{for } i = 0,\ 11\\
\ZM/3 & \text{for } i = 4,\ 8\\
\ZM/2 & \text{for } i = 6\\
0 & \text{otherwise}
\end{cases}
\]

\section{Another method for type $A$}\label{a bis}

Here we will explain a method which applies only in type $A$.
This is because the minimal class is a Richardson class only
in type $A$.

So suppose we are in type $A_{n - 1}$. We can assume $G = GL_n$.
The minimal class corresponds to the partition $(2,1^{n - 2})$.
It consists of the nilpotent matrices of rank $1$ in $\gG\lG_n$, or,
in other words, the matrices of rank $1$ and trace $0$.

Let us consider the set $E$ of pairs $([v],x) \in \PM^{n - 1} \times \gG\lG_n$
such that $\textrm{Im}(x) \incl \CM v$ (so $x$ is either zero or of
rank $1$). Together with
the natural projection, this is a vector bundle on $\PM^{n - 1}$,
corresponding to the locally free sheaf $\EC = \OC(-1)^n$
(we have one copy of the tautological bundle for each column).

There is a trace morphism $\Tr : \EC \to \OC$. Let $\FC$ be its kernel,
and let $F$ be the corresponding sub-vector bundle of $E$.
Then $F$ consists of the pairs $([v],x)$ such that $x$ is either zero
or a nilpotent matrix of rank $1$ with image $\CM v$.
The second projection gives a morphism $\pi : E \to \overline{\OC_\mini}$,
which is a resolution of singularities, with exceptional fiber
the null section. So we have an isomorphism from $F$ minus the null section
onto $\OC_\mini$.

As before, we have a Gysin exact sequence
\[
H^{i-2n+2} \elem{c} H^i \longto H^i(\OC_\mini,\ZM)
\longto H^{i-2n+3} \elem{c} H^{i+1}
\]
where $H^j$ stands for $H^j(\PM^{n-1},\ZM)$ and $c$ is the multiplication
by the last Chern class $c$ of $F$. Thus $H^i(\OC_\mini,\ZM)$ fits
in a short exact sequence
\[
0 \longto \Coker(c:H^{i-2n+2}\to H^i) \longto H^i(\OC_\mini,\ZM)
\longto \Ker(c:H^{i-2n+3}\to H^{i+1}) \longto 0
\]
We denote by $y\in H^2(\PM^{n - 1}, \ZM)$ the first Chern class of $\OC(-1)$.
We have $H^*(\PM^{n - 1},\ZM) \simeq \ZM[y]/(y^n)$ as a ring.
In particular, the cohomology of $\PM^{n-1}$ is free and
concentrated in even degrees.

For $0 \leqslant i \leqslant 2n - 4$, we have
$H^i(\OC_\mini,\ZM) \simeq H^i$ which is isomorphic to $\ZM$ if
$i$ is even, and to $0$ if $i$ is odd.
We have $H^{2n - 3}(\OC_\mini,\ZM) \simeq \Ker(c:H^0\to H^{2n-2})$
and $H^{2n - 2}(\OC_\mini,\ZM) \simeq \Coker(c:H^0\to H^{2n-2})$.
For $2n - 1 \leqslant i \leqslant 4n - 5$, we have
$H^i(\OC_\mini,\ZM) \simeq H^{i - 2n + 3}$ which is isomorphic to $\ZM$ if
$i$ is odd, and to $0$ if $i$ is even.

We have an exact sequence
\[
0 \longto \FC \longto \EC = \OC(-1)^n \longto \OC \longto 0
\]

The total Chern class of $F$ is thus
$$(1+y)^n = \sum_{i = 0}^{n-1} \binom{n}{i} y^i$$
by multiplicativity (remember that $y^n = 0$).
So its last Chern class $c$ is $ny^{n - 1}$.

In fact, $F$ can be identified with the cotangent bundle
$T^*(G/Q)$, where $Q$ is the parabolic subgroup which stabilizes
a line in $\CM^n$, and $G/Q \simeq \PM^{n - 1}$~; then
we can use the fact that the Euler characteristic of $\PM^{n - 1}$
is $n$.

We can now determine the two remaining cohomology groups.
\[
H^{2n - 3}(\OC_\mini,\ZM) \simeq \Ker(\ZM\elem{n}\ZM) = 0
\]
and
\[
H^{2n - 2}(\OC_\mini,\ZM) \simeq \Coker(\ZM\elem{n}\ZM) = \ZM/n
\]

Thus we find the same result as in section \ref{case by case} for
the cohomology of $\OC_\mini$ in type $A_{n - 1}$.

\section*{Acknowledgments}

I thank the
\'Ecole Polytechnique F\'ed\'erale de Lausanne
(EPFL)
for allowing me to attend the special semester
``Group representation theory'' in 2005. Some period was devoted to
topology and representations, and I greatly benefited from many
conversations with the participants, in particular
Alejandro Adem,
Jesper Grodal,
Hans-Werner Henn,
Ran Levi
and
Peter Symonds.
There I learned the topological tools that I needed for this work,
and I completed section 4.

Most of the writing was done at the University of Oxford, thanks to a grant
of the European Union (through the Liegrits network).

I also wish to thank C\'edric Bonnaf\'e for introducing me to the combinatorics
of section \ref{philg xi}, and Rapha\"el Rouquier for fruitful discussions.

\newcommand{\etalchar}[1]{$^{#1}$}
\def\cprime{$'$}

\noindent{\sc Daniel Juteau}\\
Institut de Math\'ematiques de Jussieu\\
Universit\'e Denis Diderot - Paris VII\\
2 place Jussieu, 75251 Paris Cedex 05, France\\
\verb|juteau@math.jussieu.fr| 

\end{document}